\documentclass[a4paper,12pt,twoside]{article}
\usepackage[francais]{babel}
\usepackage{amssymb, amsmath, eucal}
\usepackage[all]{xy}
\usepackage{mathrsfs}
\begin{document}
\textwidth15.5cm
\textheight22.5cm
\voffset=-13mm
\newtheorem{The}{Theorem}[section]
\newtheorem{Lem}[The]{Lemma}
\newtheorem{Prop}[The]{Proposition}
\newtheorem{Cor}[The]{Corollary}
\newtheorem{Rem}[The]{Remark}
\newtheorem{Titre}[The]{\!\!\!\! }
\newtheorem{Conj}[The]{Conjecture}
\newtheorem{Question}[The]{Question}
\newtheorem{Prob}[The]{Problem}
\newtheorem{Def}[The]{Definition}
\newtheorem{Not}[The]{Notation}
\newtheorem{Claim}[The]{Claim}
\newtheorem{Ex}[The]{Example}
\newcommand{\C}{\mathbb{C}}
\newcommand{\R}{\mathbb{R}}
\newcommand{\N}{\mathbb{N}}
\newcommand{\Z}{\mathbb{Z}}
\newcommand{\Proj}{\mathbb{P}}

\begin{center}

{\Large\bf Regularisation of Currents with Mass Control and Singular Morse Inequalities}

\end{center}

\begin{center}

 {\large Dan Popovici}

\end{center}

\vspace{3ex}

\noindent {\small {\bf Abstract.} Let $X$ be a compact complex, not necessarily K\"ahler, manifold of dimension $n$. We characterise the volume of any holomorphic line bundle $L\rightarrow X$ as the supremum of the Monge-Amp\`ere masses $\int_X T_{ac}^n$ over all closed positive currents $T$ in the first Chern class of $L$, where $T_{ac}$ is the absolutely continuous part of $T$ in its Lebesgue decomposition. This result, new in the non-K\"ahler context, can be seen as holomorphic Morse inequalities for the cohomology of high tensor powers of line bundles endowed with arbitrarily singular Hermitian metrics. It gives, in particular, a new bigness criterion for line bundles in terms of existence of singular Hermitian metrics satisfying positivity conditions. The proof is based on the construction of a new regularisation for closed $(1,\,1)$-currents with a control of the Monge-Amp\`ere masses of the approximating sequence. To this end, we prove a potential-theoretic result in one complex variable and study the growth of multiplier ideal sheaves associated with increasingly singular metrics.}

\vspace{3ex}

\section{Introduction}\label{subsection:introduction1}

 Let $\varphi : \Omega \rightarrow \R \cup \{-\infty\}$ be a plurisubharmonic (psh) function on an open subset $\Omega \subset \C^n$, and let $z=(z_1, \dots , \, z_n)$ be the standard coordinates on $\C^n$. The Lelong number $\nu(\varphi, \, x)$ of $\varphi$ at an arbitrary point $x\in \Omega$ is defined as the mass carried by the positive measure $dd^c\varphi \wedge (dd^c\log|z-x|)^{n-1}$ at $x$ (see, for instance, Demailly's book [Dem97], chapter III). It is a well-known result of Skoda ([Sko72a]) that the Lelong numbers of $\varphi$ affect the local integrability of $e^{-2\varphi}.$ Indeed, if $\nu(\varphi, \, x) < 1$, then $e^{-2\varphi}$ is integrable on some neighbourhood of $x$. On the contrary, if $\nu(\varphi, \, x) \geq n$, then $e^{-2\varphi}$ is not integrable near $x$. The integrability of  $e^{-2\varphi}$ is unpredictable when $1 \leq \nu(\varphi, \, x) < n$.

  Our first aim is to establish a potential-theoretic result in the case $n=1$ where there is no unpredictability interval. Let $U \subset \C$ be an open set, $\varphi_0:U \rightarrow \R \cup \{-\infty\}$ a subharmonic function, and $T=dd^c\varphi_0$ the associated closed positive current of bidegree $(1, 1)$. The current $T$ can be identified with the Laplacian $\Delta \varphi_0$ of $\varphi_0$ computed in the sense of distributions. It defines a positive measure $\mu = dd^c\varphi_0$ on $U$. In one complex variable, the mass of $dd^c\varphi_0$ at a point $x$ coincides with the Lelong number $\nu(\varphi_0, \, x)$. Let $D(x_0, \, r) \Subset U$ be an arbitrary disc of radius $0 <r<\frac{1}{2},$ and let \\

\hspace{10ex} $\displaystyle \gamma = \int\limits_{P(x_0, \, 2r)}dd^c\varphi_0$, \\

\noindent be the mass carried by the measure $dd^c\varphi_0$ on the square $P(x_0, \, 2r)$ of edge $2r$ centred at $x_0$. Consider the decomposition:

\vspace{1ex}

\hspace{6ex}     $\varphi_0 = N \star \Delta \varphi_0 + h_0,$   \hspace{2ex} on $D(x_0, \, r),$

\vspace{1ex}

\noindent where $N(z) = \frac{1}{2\pi} \, \log |z|$ is the Newton kernel in one complex variable, and $h_0 = \mathrm{Re}\, g_0$ is a harmonic function expressed as the real part of a holomorphic function $g_0.$

 First, we will be considering a method of neutralising the $-\infty$-poles of $\varphi_0$ on a fixed disc in order to make the exponential $e^{-2m\varphi_0}$ integrable and to control the growth rate of its integral as $m\rightarrow +\infty$.

\begin{The}\label{The:potentiel}

  Let $\varphi_0  : U \rightarrow \R \cup \{-\infty\}$ be a subharmonic function on an open set $U \subset \C$, and let $D\subset D(x_0, \, r) \Subset U$ be an open subset contained in a disc of radius $0 <r<\frac{1}{2}$. Fix $0<\delta <1$. Then, for every $m\gg 1$, there exist finitely many points $a_1=a_1(m), \dots , \, a_{N_m} = a_{N_m}(m)\in D,$ such that the positive integers $m_j$  defined as: \\

 $m_j=\max\{[m\nu(\varphi_0, a_j)], \,\, 1\}$, \hspace{2ex} $j= 1, \dots , N_ m,$ \hspace{2ex} ($[\hspace{1ex}]$ is the integer part), \\

\noindent  and the holomorphic function $f_m(z) = e^{m\, g_0(z)} \, \overset{N_m}{\underset{j=1}{\prod}}(z-a_j)^{m_j} $ defined on $D,$ satisfy the following properties: \\

\noindent $(i)$ \, $\overset{N_m}{\underset{j=1}{\sum}}m_j \leq m\gamma (1 + \delta)$, \hspace{2ex} where $\gamma$ is the $dd^c\varphi_0$-mass of $P(x_0, \, 2r)$ ;  \\

\noindent $(ii)$ \,  there exists a constant $C = C(r) > 0,$ independent of $m,$ such that: \\

\hspace{10ex} $\displaystyle  |a_j -  a_k| \geq \frac{C}{m^2}$,

\vspace{1ex} 

\noindent for all $a_j, \, a_k,$ such that $j\neq k$ and $\nu(\varphi_0, \, a_j),  \nu(\varphi_0, \, a_k) < \frac{1-\delta}{m}; $ \\

\noindent $(iii)$ \, $\displaystyle \int_{D}|f_m(z)|^2 e^{-2m\varphi_0(z)}\, d\lambda(z)=o(m), \hspace{2ex} \mathrm{when} \hspace{2ex} m\rightarrow +\infty, $ \\

\noindent where $d\lambda$ is the Lebesgue measure in $\C.$

\end{The}

 Higher dimensional analogues of this result have yet to be found. However, the Ohsawa-Takegoshi $L^2$ extension theorem (see [OT87], [Ohs88]) applied on a complex line enables us to derive geometric applications of Theorem \ref{The:potentiel} in several complex variables. The first application is a global regularisation theorem for closed almost positive $(1, \, 1)$-currents in the spirit of Demailly (see [Dem92]), but with an additional control on the Monge-Amp\`ere masses of the regularising currents. Here is the set-up. 

  Let $T$ be a $d$-closed current of bidegree $(1, \, 1)$ on  a compact complex manifold $X$ of dimension $n$. Assume that $T\geq \gamma$ for some real continuous $(1, \, 1)$-form $\gamma$ (i. e. $T$ is almost positive). The current $T$ can be globally written as $T=\alpha + dd^c\varphi$, with a global $C^{\infty}$ $(1, \, 1)$-form $\alpha$ and an almost psh potential $\varphi$ on $X$ (i.e. $\varphi$ can be locally expressed as the sum of a psh function and a $C^{\infty}$ function). The notation $dd^c:=\frac{i}{\pi}\, \partial\bar{\partial}$ will be used in all that follows. A variant of Demailly's regularisation theorem (see [Dem92, Proposition 3.7]) asserts that $T$ is the weak limit of currents $T_m=\alpha + dd^c\varphi_m$ lying in the $\partial\bar{\partial}$-cohomology class of $T$ and having analytic singularities. These are, by definition, singularities for which $\varphi_m$ can be locally written as 

\begin{equation}\label{equation:analsing}\frac{c}{2}\, \log(|g_1|^2 + \dots + |g_N|^2)+ C^{\infty},\end{equation}

\noindent with a constant $c>0$, and holomorphic functions $g_1, \dots , \, g_N$. Each $T_m$ can be chosen to be smooth on $X\setminus V{\mathcal I}(mT)$, where $V{\mathcal I}(mT)$ is the zero variety of the multiplier ideal sheaf ${\mathcal I}(mT)$ associated with $mT$ (defined locally as ${\mathcal I}(m\varphi)$, see (\ref{eqn:misdef})). Moreover, for any Hermitian metric $\omega$ on $X$, $T_m$ can be chosen such that: \\

\hspace{17ex}  $T_m \geq \gamma - \varepsilon_m \, \omega$, \hspace{2ex} for some sequence $\varepsilon_m \downarrow 0$, \\

\noindent and the Lelong numbers satisfy: $\displaystyle \nu(T, \, x)-\frac{n}{m} \leq \nu(T_m, \, x) \leq \nu(T, \, x), \hspace{2ex} x\in X.$ \\

\noindent What this theorem does not specify, however, is whether there exist regularisations $T_m\rightarrow T$ with analytic singularities having the extra property that the growth in $m$ of the masses of the wedge-power currents $T_m^k$ (Monge-Amp\`ere currents) is under control. In other words, we would like to control the growth rate of the quantities:  \\

\hspace{10ex} $\displaystyle \int\limits_{X\setminus V{\mathcal I}(mT)}(T_m-\gamma + \varepsilon_m \, \omega)^k \wedge \omega^{n-k}$, \hspace{2ex} $k=1, \dots , \, n,$ \\

\noindent as $m\rightarrow +\infty$, where $V{\mathcal I}(mT)=\{\varphi_m=-\infty\}$ is the polar set of $T_m$. Using Theorem \ref{The:potentiel} we can modify Demailly's original construction to settle this question in the following form.

\begin{The}\label{The:globalregularisation}  Let $T\geq \gamma$ be a $d$-closed current of bidegree $(1, \, 1)$ on a compact complex manifold $X$, where $\gamma$ is a continuous $(1, \, 1)$-form such that $d\gamma=0$. Then, in the $\partial\bar{\partial}$-cohomology class of $T$, there exist closed $(1, \, 1)$-currents $T_m$ with analytic singularities converging to $T$ in the weak topology of currents such that each $T_m$ is smooth on $X\setminus V{\mathcal I}(mT)$ and: \\

\noindent $(a)$ \, $\displaystyle T_m \geq \gamma - \frac{C}{m}\, \omega$, \hspace{2ex} $m\in \N$; \\

\noindent $(b)$ \, $\displaystyle \nu(T, \, x)-\varepsilon_m \leq \nu(T_m, \, x) \leq \nu(T, \, x), \hspace{2ex} x\in X, \, m\in \N$, for some $\varepsilon_m\downarrow 0$; \\

\noindent $(c)$ \,  $\displaystyle \lim\limits_{m\rightarrow +\infty}\frac{1}{m}\, \int\limits_{X\setminus V{\mathcal I}(mT)}(T_m -\gamma + \frac{C}{m}\, \omega)^k \wedge \omega^{n-k} = 0$, \hspace{2ex} $k=1, \dots , \, n=\mbox{dim}_{\C}X$, \\

\noindent where $\omega$ is an arbitrary Hermitian metric on $X$.

\end{The}




\vspace{2ex}

 This result can be used to prove a new characterisation of big line bundles in terms of curvature currents. Let us briefly review a few basic facts. A holomorphic line bundle $L$ over a compact complex manifold $X$ of dimension $n$ is said to be big if $\mbox{dim}_{\C}\, H^0(X, \, L^m) \geq C\, m^n$ for some constant $C>0$ and for all large enough $m\in \N$. This amounts to the global sections of $L^m$ defining a bimeromorphic embedding of $X$ into a projective space for $m\gg 0$.  The compact manifold $X$ is said to be Moishezon if the transcendence degree of its meromorphic function field equals $n:=\mbox{dim}_{\C}X$, or equivalently, if there exist $n$ global meromorphic functions that are algebraically independent. A Moishezon manifold becomes projective after finitely many blow-ups with smooth centres. There is, moreover, a bimeromorphic counterpart to Kodaira's embedding theorem: a compact complex manifold $X$ is Moishezon if and only if there exists a big line bundle $L\rightarrow X$. The asymptotic growth of the dimension of $H^0(X, \, L^m)$ as $m\rightarrow +\infty$ is actually measured by a birational invariant of $L$, the {\it volume}, defined as: \\

\hspace{6ex} $\displaystyle v(L):=\limsup_{m\rightarrow +\infty} \frac{n!}{m^n}\, h^0(X, \, L^m).$ \\

\noindent Clearly, $L$ is big if and only if $v(L)>0$. Switching now to the analytic point of view, recall that a singular Hermitian metric $h$ on $L$ is defined in a local trivialisation $L_{|U}\simeq U\times \C$ as $h=e^{-\varphi}$ for some weight function $\varphi :U \rightarrow [-\infty , \, +\infty )$ which is only assumed to be locally integrable. In particular, the singularity set $\{x\in U, \, \varphi(x)=-\infty\}$ is Lebesgue negligible. The associated curvature current $T:=i\Theta_h(L)$ is a closed current of bidegree $(1, \, 1)$ on $X$ representing the first Chern class $c_1(L)$ of $L$. It is locally defined as $T=dd^c\varphi$, where $\varphi$ is a local weight function of $h$. 

Recall that an almost positive current $T$ can be locally written in coordinates as $T=\sum\limits_{j, \, k} T_{j, \, k}\, dz_j \wedge d\bar{z}_k$ for some complex measures $T_{j, \, k}$. The Lebesgue decomposition of the coefficients $T_{j, \, k}$ into an absolutely continous part and a singular part with respect to the Lebesgue measure induces a current decomposition as $T=T_{ac} + T_{sing}$. By the Radon-Nikodym theorem, the coefficients of the absolutely continous part are $L^1_{loc}$ and thus the exterior powers $T_{ac}^m$ are well defined (though not necessarily closed) currents for $m=1, \dots, \, n$.

When applied to curvature currents, the current regularisation theorem \ref{The:globalregularisation} with controlled Monge-Amp\`ere masses enables us to characterise the volume of a line bundle in terms of positive currents in $c_1(L)$. This gives in particular a bigness criterion for line bundles in terms of existence of singular Hermitian metrics satisfying positivity assumptions (and implicitly a characterisation of Moishezon manifolds).

\begin{The}\label{The:bigness} Let $L$ be a holomorphic line bundle over a compact complex manifold $X$. Then the volume of $L$ is characterised as: \\

\hspace{6ex} $\displaystyle v(L)=\sup_{T\in c_1(L), \, T\geq 0} \int_X T_{ac}^n$. \\

\noindent In particular, $L$ is big if and only if there exists a possibly singular Hermitian metric $h$ on $L$ whose curvature current $T:=i\Theta_h(L)$ satisfies the following positivity conditions: 

\vspace{3ex}

 \hspace{6ex} $(i)$ \, $T\geq 0$ \hspace{2ex} on $X$; \hspace{3ex} $(ii)$ \, $\displaystyle \int_X T_{ac}^n>0$. \\

\end{The}

 In the special case when the ambient manifold $X$ is K\"ahler, this same result was obtained by Boucksom ([Bou02, Theorem 1.2]). This strengthens a previous (only sufficient) bigness criterion by Siu ([Siu85]) that solved affirmatively the Grauert-Riemenschneider conjecture ([GR70]). Bigness was guaranteed there under the extra assumption that the curvature current $T$ (or the metric $h$) be $C^{\infty}$. Theorem \ref{The:bigness} falls into the mould of ideas originating in Demailly's holomorphic Morse inequalities ([Dem85]). Its proof hinges on the regularisation Theorem \ref{The:globalregularisation} above and on Bonavero's singular version of Demailly's Morse inequalities ([Bon98]). It strengthens a bigness criterion in [Bon98] which required the curvature current $T$ to have analytic singularities. On the other hand, Ji and Shiffman ([JS93]) proved that $L$ being big is equivalent to $L$ having a singular metric whose curvature current is strictly positive on $X$, (i.e. \,\, $\geq \varepsilon \omega$ for some small $\varepsilon >0$). This implies, in particular, the ``only if'' part of the above Theorem \ref{The:bigness}. The thrust of the new ``if'' part of Theorem \ref{The:bigness} is to relax the strict positivity assumption on the curvature current.

 Let us finally stress that the main interest of Theorems \ref{The:globalregularisation} and \ref{The:bigness} lies in $X$ being an arbitrary compact manifold. Related results are known to exist for K\"ahler manifolds (e.g. [DP04, Theorem 0.4], [Bou02, Theorem 1.2]). The approach to the non-K\"ahler case treated here is quite different. The crux of the argument is modifying the existing procedure for regularising $(1, \, 1)$-currents to get an effective control on the Monge-Amp\`ere masses (Theorem \ref{The:globalregularisation}). If $X$ is K\"ahler, the sequence of masses in the usual Demailly regularisation of currents is easily seen to be bounded by applying Stokes's theorem and using the closedness of $\omega$ (see [Bou02]). The situation is vastly different in the non-K\"ahler case where a new regularisation of currents is needed with a possibly unbounded sequence of masses.

\vspace{2ex}
 
The paper runs as follows. Section \ref{section:preliminaries} collects the main ingredients for the proof of Theorem \ref{The:potentiel}. The most prominent of these is an atomisation lemma for positive measures in $\R^2$ due to Yulmukhametov ([Yul85]) and Drasin ([Dra01]). Section \ref{subsection:proofpotentiel} goes on to complete the proof of Theorem \ref{The:potentiel}.
 

After a preliminary section, the next three sections achieve the proof of Theorem \ref{The:globalregularisation}. Section \ref{subsection:modreg} explains the new regularisation technique for psh functions on domains in $\C^n$ based on using derivatives of holomorphic functions in a weighted Bergman space. Section \ref{subsection:addefect} introduces another key idea of this paper, namely an effective control, with estimates, of how far multiplier ideal sheaves are from behaving linearly as the singularities increase. It turns out that in the case of analytic singularities the growth of multiplier ideal sheaves is almost linear. This provides a useful complement to the Demailly-Ein-Lazarsfeld subadditivity theorem. Section \ref{subsection:reg} combines ideas of the previous sections to complete the proof of Theorem \ref{The:globalregularisation}. 

Finally, section \ref{section:bigness} gives a proof of Theorem \ref{The:bigness} as a geometric application of the new regularisation procedure for currents. An appendix, section \ref{subsection:appendix}, provides a few explanations that have fallen between the cracks of the previous sections. 

\vspace{3ex}

\noindent {\bf Acknowledgments.} The author is extremely grateful to Jean-Pierre Demailly for suggesting the problem and for numerous helpful discussions on the topic. The author would like to thank the referee for pertinent remarks and comments. A word of thanks is also due to quite a number of people who showed an interest in this work, of whom Bo Berndtsson, S\'ebastien Boucksom and Takeo Ohsawa are only a few.

\section{Preliminaries to Theorem \ref{The:potentiel}}\label{section:preliminaries}

 In this section we clear the way to the proof of Theorem \ref{The:potentiel}. The set-up is the one described in the Introduction. Fix $m\in \mathbb{N}^{\star}$ and $\delta>0$. As the upper level set for Lelong numbers:\\

 \hspace{6ex} $E_{1-\delta}(m\,dd^c\varphi_0)\! :=\{ x\in U \, ; \,  \nu(m\,\varphi_0, \, x)\geq 1-\delta \}$ \\

\noindent is analytic of dimension $0$, its intersection with a relatively compact subset is finite. Let  \\

\hspace{6ex} $E_{1-\delta}(m\,dd^c\varphi_0) \cap D: =\{ a_1, \dots , \, a_{p(m)} \}. $ \\

\noindent As $dd^c\psi_m = m \, dd^c\varphi_0 - \overset{p(m)}{\underset{j=1}\sum}[m \, \nu(\varphi_0, \, a_j)] \, \delta_{a_j}  \geq 0$ as currents, the function

\vspace{1ex}

\hspace{6ex} $\psi_m(z):=m \, \varphi_0(z) - \overset{p(m)}{\underset{j=1}\sum} [m \, \nu(\varphi_0, \, a_j)] \, \log|z-a_j|$

\vspace{1ex}

\noindent is still subharmonic and $\nu(\psi_m, \, a_j)= m \, \nu(\varphi_0, \, a_j) -[m \, \nu(\varphi_0, \, a_j)]$ for every $j$. \\

\noindent To ensure the integrability property $(iii)$ in Theorem \ref{The:potentiel}, we must introduce a factor $(z-a_j)^{m_j}$ in the definition of the function $f_m$ being constructed for every point $a_j$ where $m_j:=[m \, \nu(\varphi_0, \, a_j)] \geq 1$. To find the other points, we can assume in all that follows, after replacing $m\varphi_0$ with $\psi_m$, that:

\begin{equation}\label{eqn:hyp}m\, \nu(\varphi_0, \, x) < 1, \hspace{2ex} \mbox{for all} \,\, x\in D.\end{equation}

\noindent Once the point masses of the current $m \, dd^c\varphi_0$ have been brought down below $1$, we still have to neutralise any diffuse mass scattered over the domain $D$ which could prevent the integral of $|f_m|^2 \, e^{-2m\varphi_0}$ from having the desired slow growth in $m$. The following lemma gives an upper bound for $e^{-2\varphi_0}$ in terms of the mass of the associated current $dd^c\varphi_0$.

\begin{Lem}\label{Lem:integetmasse} With the notation in the introduction, if $\gamma: = \int_{D}dd^c\varphi_0$, the following estimate holds: \\

\hspace{6ex} $\displaystyle e^{-2(\varphi_0(z) - h_0(z))} \leq \frac{1}{\int\limits_ D dd^c \varphi_0} \, \int\limits_ D \frac{1}{|\zeta-z|^{2\gamma}} \, dd^c\varphi_0(\zeta),$

\noindent for all $z\in D$.

\end{Lem}

\noindent {\it Proof.} Let $d\mu(\zeta): = \gamma^{-1} \, dd^c\varphi_0(\zeta)$ be a probability measure on $D.$ For all $z \in D,$ we have: \\

$(\varphi_0 - h_0)(z) = \displaystyle \int_{D}\log|\zeta - z| \, dd^c\varphi_0(\zeta),$ \\

\noindent or, equivalently, \\

$-(\varphi_0 - h_0)(z) = \displaystyle \int_{D} \gamma  \log|\zeta - z|^{-1} d\mu(\zeta),$  \hspace{2ex} $z\in D.$ \\

\noindent Now, Jensen's convexity inequality entails: \\

$\displaystyle e^{-2(\varphi_0 - h_0)(z)} \leq \int_{D} e^{2\gamma  \log|\zeta - z|^{-1}} d\mu(\zeta)= \gamma^{-1} \int_{D}\frac{1}{|\zeta-z|^{2\gamma}}dd^c\varphi_0(\zeta),$ \\

\noindent which proves the lemma.   \hfill $\Box$

\vspace{3ex}

\noindent Applying Lemma \ref{Lem:integetmasse} to the function $m \, \varphi_0,$ we get: \\

 $\displaystyle e^{-2m(\varphi_0(z)-h_0(z))}\leq \frac{1}{\int\limits_{D}dd^c\varphi_0}\, \int\limits_{D}\frac{1}{|\zeta - z|^{2m\gamma}} \, dd^c\varphi_0(\zeta),$ \hspace{2ex} $z\in D$. \\

\noindent The right-hand term above may not be integrable as a function of $z$ when $m\gamma > 1$. To get around this, we will cut $D$ into pieces in such a way that each piece has a mass $< 1$ for the measure $m \, dd^c\varphi_0$ and the number of pieces does not exceed $m\gamma (1+\delta)$. We will subsequently choose a point in each piece, intuitively its ``centre'', and will define $f_m$ as a holomorphic function on $D$ whose only zeroes are these points. This will be shown to satisfy the conditions in Theorem \ref{The:potentiel}. 

 Cutting $D$ into pieces relies on the following lemma due to Yulmukhametov ([Yul85]) and, in a generalised form, to Drasin ([Dra01, Theorem 2.1]). It describes an atomisation procedure for arbitrary positive measures $\mu$ in one complex variable. This is the main technical ingredient in the proof of Theorem \ref{The:potentiel}.

\begin{Lem}{\it ([Yul85]; [Dra01].)}\label{Lem:atomisation} Let $\mu$ be a positive measure supported in a square $R \subset \R^2$ with sides parallel to the coordinate axes. Suppose $\mu(R)= N >1, \, N\in \mathbb{Z}.$ Then, there exist a family of closed rectangles $(R_j)_{1\leq j \leq N}$ with sides parallel to the coordinate axes, and a family of positive measures $(\mu_j)_{1\leq j \leq N},$ such that:  \\

\noindent (a) \, $\mu = \overset{N}{\underset{j=1}{\sum}}\mu_j$, \, $\mu_j(\R^2) = 1$ and \, $\mathrm{Supp}\,\mu_j \subset R_j$ \,\, for all \,\, $j=1, \dots , N;$  \\

\noindent (b) \, $R = \overset{N}{\underset{j=1}{\bigcup}}R_j = \overset{N}{\underset{j=1}{\bigcup}}\mathrm{Supp}\,\mu_j;$ \\

\noindent (c) \, the interiors of the convex hulls of the supports $\mathrm{Supp}\,\mu_j$ of the $\mu_j$'s are mutually disjoint;  \\

\noindent (d) \, the ratio of the sides of each rectangle $R_j$ lies in the interval $[\frac{1}{3}, \, 3]$ (i. e. $R_j$ is an ``almost square'' in the terminology of [Dra01]); \\

\noindent (e) \, each point in $\R^2$ belongs to the interior of at most four distinct rectangles $R_j;$  \\

\noindent (f) \, each $\mathrm{Supp}\,\mu_j$ is contained in a rectangle $P_j\subset R_j$, and the distance between the centres of any two distinct rectangles $P_j$ is $\geq \frac{C}{N^2},$ where $C>0$ is the side of the square $R$.

\end{Lem}

\noindent {\it Idea of proof (according to [Dra01]).} Yulmukhametov originally proved this result (see [Yul85]) for absolutely continuous measures $\mu$. The generalisation to the case of arbitrary measures is due to Drasin ([Dra01]). We summarise here the ideas of Drasin's proof. Conclusion $(f)$ was not explicitly stated, but it can be easily inferred from the proof given there. The first idea is to reduce the problem to the case of a measure $\mu$ satisfying $\mu(p) < 1$ at every point $p\in R.$ This is done by subtracting from the original measure $\mu$ the integer part $[\mu(p)]$ of each point mass $\mu(p) > 1.$ We may also assume, after a possible rotation of the coordinate system of $\R^2,$ that for every line $L$ parallel to one of the coordinate axes, there exists at most one point $p\in L$ such that $\mu(p) > 0$ while $\mu(L \setminus p) = 0.$ 

  After these reductions, the key step is to prove that if an almost square $R$ contains the support of a measure $\mu$ satisfying these properties, then there exist almost squares $R_0$ and $R_1$ and a decomposition $\mu=\mu_0 + \mu_1$ such that $\mathrm{Supp}\, \mu_j \subset R_j,$ $j=0, \, 1,$ which satisfies conclusions $(b) - (d)$ of the lemma. The masses $\mu_j(R_j)$ are integers. If $\mu_j(R_j) > 1,$ we repeat this procedure to obtain almost squares $R_{j, \, 0}$, $R_{j, \, 1}$ and a decomposition $\mu_j = \mu_{j, \, 0} + \mu_{j, \, 1}.$ By repeatedly applying this procedure we get almost squares $R_I$ and measures $\mu_I,$ indexed over multi-indices $I=i_1, \dots, \, i_l$ made up of digits $0$ and $1$. The procedure terminates when all masses $\mu_I(R_I)=1.$ A technical lemma then yields conclusion $(e)$ and thus clinches the proof of this result. We refer for details to Drasin ([Dra01, $\S 2$, p. 165-171]).  \hfill $\Box$

\vspace{2ex}

\section {Proof of Theorem \ref{The:potentiel}.}\label{subsection:proofpotentiel} Building on preliminaries in the previous section, we will now complete the proof of Theorem \ref{The:potentiel}. The notation and set-up are unchanged. Let $R$ be the square of edge $2r$ centered at $x_0$. It contains $D(x_0, \, r)$ and implicitly $D$. By hypothesis (\ref{eqn:hyp}), we may assume that $m\, \nu(\varphi_0, \, x) < 1$ for all $x\in D$. The positive measure $\mu:= dd^c\varphi_0$ has mass $\gamma$ on $R$. Fix $0< \delta < 1$ and, for $m \gg  1,$ choose an integer $N_m$ such that:

\begin{equation}\label{eqn:choiceofno}\frac{2}{2-\delta}\,  m\, \gamma < N_m \leq m\, \gamma \, (1+\delta).\end{equation}

\noindent Such an integer exists if $m$ is so large that $m \, \gamma \, (1+\delta - \frac{2}{2-\delta}) = m\, \gamma \, \frac{\delta(1-\delta)}{2-\delta} > 1$. We now apply the atomisation Lemma \ref{Lem:atomisation} to the measure $\frac{N_m}{\gamma}\, \mu = \frac{N_m}{\gamma}\, dd^c\varphi_0$ of total mass $N:=N_m$ on $R$. We get a covering of $R$ by closed rectangles: \\

\hspace{6ex} $R=\overset{N_m}{\underset{j=1}{\bigcup}}R_j(m),$ \\

\noindent and a decomposition of measures $\displaystyle\frac{N_m}{\gamma}\, \mu = \overset{N_m}{\underset{j=1}\sum} \nu_{m, \, j}$ such that every $R_j:=R_j(m)$ is an almost square containing the support of $\nu_{m, \, j}$, and $\nu_{m, \, j}(R_j)=1.$ \\

\noindent By part $(f)$ of Lemma \ref{Lem:atomisation}, there is a family of possibly smaller rectangles $P_j(m)\subset R_j(m)$, $j=1, \dots , N_m,$ still covering $R$ such that each $P_j(m)$ contains the support of $\nu_{m, \, j}$. Their main feature is the following: if $a_j=a_j(m)$ is the centre of $P_j(m)$, the mutual distances of the $a_j$'s are $\geq C/N_m^2$ with $C>0$ independent of $m$. Unlike the $R_j(m)$'s, we do not know whether the $P_j(m)$'s are almost squares. \\

We will prove that the integer $N_m$ and the points $a_j$ satisfy the conclusions of Theorem \ref{The:potentiel}. As by Hypothesis (\ref{eqn:hyp}) we have: \\

\hspace{6ex} $m_j= \max\{[m \, \nu(\varphi_0, \, a_j)], \, 1\} = 1,$ \hspace{2ex} for $j=1, \dots , \, N_m$,

\vspace{2ex}

\noindent we get: $\overset{N_m}{\underset{j=1}{\sum}} m_j = N_m \leq m  \gamma  (1 + \delta),$ which is conclusion $(i)$ of Theorem  \ref{The:potentiel}. The lower bound on the mutual distances of the points $a_j$ and the choice of $N_m=O(m)$ (cf. (\ref{eqn:choiceofno})) ensure that the points $a_j$ satisfy the conclusion $(ii)$ of Theorem \ref{The:potentiel}. It remains to check that for the holomorphic function:

\vspace{2ex}

 \hspace{6ex} $f_m(z):= e^{m\, g_0(z)}\, \overset{N_m}{\underset{j=1}{\prod}}(z-a_j),$  \hspace{2ex}  $z\in D,$

\vspace{2ex}

\noindent the integral $\int\limits_{D}|f_m|^2 \,  e^{-2m\varphi_0}\, d\lambda$ has the desired slow growth in $m$ as $m\rightarrow +\infty$. Since

\vspace{1ex}

\hspace{6ex} $\displaystyle \int\limits_{D}|f_m|^2 \, e^{-2m\varphi_0}\, d\lambda \leq \overset{N_m}{\underset{j=1}{\sum}}\int\limits_{R_j}|f_m|^2 \, e^{-2m\varphi_0} \, d\lambda,$

\vspace{1ex}

\noindent the analysis is reduced to finding a convenient upper bound for each integral on $R_j$. Fix $j\in \{1, \dots ,N_m\}$. Since $a_j\in R_j$, we see that $R_j$ is contained in the disc $D_j:=D(a_j, \, r_j)$, where $r_j$ is $\sqrt{2}$ times the longest edge of $R_j$. Conclusion $(e)$ of Lemma \ref{Lem:atomisation} implies that the sum of the Euclidian areas of the $R_j$'s is bounded above by four times the area of the square $R$ of edge $2r$. As the $R_j$'s are almost squares, it follows that there is a constant $C_1(r)>0,$ depending only on $r$, such that

\vspace{1ex}

 \hspace{10ex}  $(e')$  \hspace{2ex}  $\overset{N_m}{\underset{j=1}\sum} r_j^2 \leq C_1(r),$  \hspace{2ex} for all $m\gg 1.$

\vspace{1ex}

\noindent Lemma \ref{Lem:integetmasse}, when applied to the function $m \, \varphi_0$ on $D_j=D(a_j, \, r_j)$, gives:

\vspace{2ex}

\noindent $\begin{array}{lll}

|f_m(z)|^2 \, e^{-2m\varphi_0(z)} & = & \overset{N_m}{\underset{k=1}\prod} |z-a_k|^2 \, e^{-2m(\varphi_0(z) - h_0(z))} \\

\vspace{2ex}

   & \leq & \displaystyle (2r)^{2(N_m-1)} \, |z-a_j|^2 \, \frac{1}{\int\limits_{D_j}dd^c\varphi_0} \,\int\limits_{D_j}\frac{1}{|\zeta - z|^{2\frac{m\gamma}{N_m}}} \,  dd^c\varphi_0(\zeta),

\end{array}$

\noindent for all $ z\in D_j.$ We have used the obvious upper bound $|z-a_k|^2 \leq (2r)^2$, for all $k\neq j.$ When integrating above with respect to $z \in D_j,$ Fubini's theorem yields: 

\begin{equation}\label{eqn:Fubini}\displaystyle \int\limits_{D_j}|f_m(z)|^2 \, e^{-2m\varphi_0(z)} d\lambda(z) \leq \frac{(2r)^{2(N_m-1)}}{\int\limits_{D_j}dd^c\varphi_0}  \int\limits_{D_j}\bigg( \int\limits_{D_j}\frac{|z-a_j|^2}{|z-\zeta|^{2\frac{m\gamma}{N_m}}}d\lambda(z) \bigg) dd^c\varphi_0(\zeta).\end{equation}

\noindent  Let us now concentrate on the integral in $z$ on the right-hand side. We get the following estimate for every $\zeta \in D_j:$

 \noindent\begin{eqnarray}\nonumber & \displaystyle \int\limits_{D_j} & \frac{|z-a_j|^2}{|z-\zeta|^{2\frac{m\gamma}{N_m}}}d\lambda(z)
 =  \int\limits_{D(a_j, \, r_j)}\frac{|z-a_j|^2}{|(z-a_j)-(\zeta-a_j)|^{2\frac{m\gamma}{N_m}}} d\lambda (z-a_j)\\
\label{eqn:rhs} & \leq &  4\pi (|\zeta-a_j|+r_j)^{2(1-\frac{m\gamma}{N_m})} \, \left(\frac{(|\zeta - a_j|+r_j)^2}{2(2-\frac{m\gamma}{N_m})} + \frac{|\zeta - a_j|^2}{2(1-\frac{m\gamma}{N_m})}\right).
\end{eqnarray}

\noindent Indeed, if we make the change of variable $x=z-a_j$ and set $\zeta -a_j=a$, we are reduced to estimating the integral: \\

\hspace{15ex}$\displaystyle\int\limits_{D(0, \, r)}\frac{|x|^2}{|x-a|^{\tau}}d\lambda(x),$

\noindent where we have set $r_j:=r$ and $\tau:=2\frac{m\gamma}{N_m}$ to simplify the notation. By the choice (\ref{eqn:choiceofno}) of $N_m,$ we have: $0 < \tau < 2.$ The change of variable $x-a=y$, followed by a switch to polar coordinates with $|y|=\rho$, implies:

 \begin{eqnarray*}\int\limits_{D(0, \, r)}\frac{|x|^2}{|x-a|^{\tau}}d\lambda(x) & = & \int\limits_{D(-a, \, r)}\frac{|y+a|^2}{|y|^{\tau}} d\lambda(y) \leq \int\limits_{D(-a, \, r)}\frac{(|y|+|a|)^2}{|y|^{\tau}}d\lambda(y)  \\
 & \leq & 2 \int\limits_{D(-a, \, r)}\frac{|y|^2+|a|^2}{|y|^{\tau}}d\lambda(y) 
\end{eqnarray*}

 \begin{eqnarray*} & = & 2\int\limits_{D(-a, \, r)}|y|^{2-\tau}d\lambda(y) + 2|a|^2\int\limits_{D(-a, \, r)}|y|^{-\tau}d\lambda(y)  \\
& \leq & 2\pi\, \left( 2\int_0^{|a| +r}\rho^{2-\tau}\rho\, d\rho + 2|a|^2\int_0^{|a|+r}\rho^{-\tau}\rho\, d\rho \right) \\
& = & 4\pi \left( |a|+r \right)^{2-\tau}\left( \frac{(|a|+r)^2}{4-\tau}+ \frac{|a|^2}{2-\tau} \right).
\end{eqnarray*}


\noindent For $r=r_j$, this gives the estimate (\ref{eqn:rhs}). Now relations (\ref{eqn:Fubini}) and (\ref{eqn:rhs}) imply: \\

\noindent$\begin{array}{lll} 

\displaystyle \int\limits_{D_j}|f_m(z)|^2 \, e^{-2m\varphi_0(z)} d\lambda(z) \leq \frac{4\pi}{\int_{D_j}dd^c\varphi_0}(2r)^{2(N_m-1)}\cdot &  &  \\
 
\hspace{6ex} \displaystyle \cdot \int_{D_j}(|\zeta - a_j|+ r_j )^{2(1-\frac{m\gamma}{N_m})} \left( \frac{(|\zeta - a_j| + r_j)^2}{2(2-\frac{m\gamma}{N_m})} + \frac{|\zeta - a_j|^2}{2(1-\frac{m\gamma}{N_m})} \right) dd^c\varphi_0(\zeta). &  &   

\end{array}$

\vspace{2ex}

\noindent Let us now shift to polar coordinates with $|\zeta - a_j|=\rho$. This gives $dd^c\varphi_0(\zeta)=dn(\rho)$, where $n(\rho)=\int_{D(a_j, \, \rho)}dd^c\varphi_0$, for all $\rho \geq 0$. Since $D_j$ is assumed to be $D(a_j, \, r_j)$, we get: \\

$\begin{array}{lll} 

\displaystyle \int\limits_{D_j}|f_m(z)|^2 \, e^{-2m\varphi_0(z)} d\lambda(z) \leq & & \\

\displaystyle \leq C(r, \, r_j) \, \int_0^{r_j}(\rho + r_j)^{2(1-\frac{m\gamma}{N_m})} \left( \frac{(\rho + r_j)^2}{2(2-\frac{m\gamma}{N_m})}+ \frac{\rho^2}{2(1-\frac{m\gamma}{N_m})} \right)n'(\rho)\, d\rho, & &

\end{array}$

\vspace{2ex}

\noindent where $\displaystyle C(r, \, r_j)= \frac{8\pi^2}{\int_{D_j}dd^c\varphi_0}(2r)^{2(N_m-1)}$. The last expression can be successively written as:

\renewcommand{\arraystretch}{3}

\noindent$\begin{array}{lll}  

\displaystyle \frac{C(r, \, r_j)}{2(2-\frac{m\gamma}{N_m})} \, \int_0^{r_j}(\rho + r_j)^{2(2-\frac{m\gamma}{N_m})}n'(\rho)\, d\rho + \frac{C(r, \, r_j)}{2(1- \frac{m\gamma}{N_m})} \, \int_0^{r_j}\rho^2(\rho + r_j)^{2(1-\frac{m\gamma}{N_m})}n'(\rho) \, d\rho  &  &  \\ 

\displaystyle = \frac{C(r, \, r_j)}{2(2-\frac{m\gamma}{N_m})} \, \bigg( n(r_j)(2r_j)^{2(2-\frac{m\gamma}{N_m})} - 2 \left( 2-\frac{m\gamma}{N_m} \right) \int_0^{r_j}n(\rho)(\rho + r_j)^{3-2\frac{m\gamma}{N_m}}\, d\rho \bigg) + &  &  \\ 

\displaystyle + \frac{C(r, \, r_j)}{2(1- \frac{m\gamma}{N_m})} \, \bigg( n(r_j)r_j^2 (2r_j)^{2(1-\frac{m\gamma}{N_m})} + &  &  \\

\displaystyle \hspace{3em} - \int_0^{r_j}n(\rho)\left[ 2\rho(\rho + r_j)^{2(1-\frac{m\gamma}{N_m})} + 2 \left( 1-\frac{m\gamma}{N_m} \right) (\rho + r_j)^{1-2\frac{m\gamma}{N_m}} \rho^2 \right] \, d\rho \bigg).  &  & 

\end{array}$

\vspace{2ex}

\noindent The terms appearing above with a `` $-$'' sign are all negative since $1-\frac{m\gamma}{N_m} > 0$ (and implicitly $2-\frac{m\gamma}{N_m} > 0$). Therefore, they can be ignored. We thus get the following upper estimate: \\

\renewcommand{\arraystretch}{1.3}

\noindent$\begin{array}{lll}

\displaystyle \int\limits_{D_j}|f_m(z)|^2 \, & e^{-2m\varphi_0(z)}  d\lambda(z)  \leq &   \\

\displaystyle  & \leq  C(r, \, r_j) \, n(r_j) \, \displaystyle \bigg( \frac{(2r_j)^{2(2-\frac{m\gamma}{N_m})}}{2(2-\frac{m\gamma}{N_m})} + \frac{r_j^2(2r_j)^{2(1-\frac{m\gamma}{N_m})}}{2(1-\frac{m\gamma}{N_m})} \bigg). &

\end{array}$

\vspace{2ex}

\noindent Since $n(r_j)=\int_{D_j}dd^c\varphi_0$, the previous upper bound and the formula of $C(r, \, r_j)$ show that:

\begin{equation}\label{eqn:pub}\int\limits_{D_j}|f_m(z)|^2 \, e^{-2m\varphi_0(z)} d\lambda(z) \leq  C\left(r, \, \frac{m\gamma}{N_m}\right) \cdot r_j^{2(2-\frac{m\gamma}{N_m})},\end{equation}

\noindent where the constant $C(r, \, \frac{m\gamma}{N_m})$ is given by the formula:

\vspace{2ex}

$$ C\left(r, \, \frac{m\gamma}{N_m}\right)=8\pi^2 \, (2r)^{2(N_m-1)}\, \bigg(\frac{2^{2(2-\frac{m\gamma}{N_m})}}{2(2-\frac{m\gamma}{N_m})} + \frac{2^{2(1-\frac{m\gamma}{N_m})}}{2(1-\frac{m\gamma}{N_m})} \bigg).  $$

\vspace{2ex}

\noindent Since estimate (\ref{eqn:pub}) holds for all indices $j\in \{1, \dots , \, N_m \}$, we get, after summing over $j$, that:

$$\int_{D}|f_m(z)|^2  \, e^{-2m\varphi_0(z)} \, d\lambda (z) \leq C\big(r, \, \frac{m\gamma}{N_m}\big) \, \overset{N_m}{\underset{j=1}{\sum}}r_j^{2(2-\frac{m\gamma}{N_m})}. $$

\vspace{2ex}

\noindent  The choice of $N_m$ was made in such a way that $1-\delta < \frac{1}{1+\delta} \leq \frac{m\gamma}{N_m} < 1 - \frac{\delta}{2}$ (cf. (\ref{eqn:choiceofno})),  which implies:

\vspace{1ex}

 \hspace{6ex}$\displaystyle \frac{\delta}{2} <  1-\frac{m\gamma}{N_m} < \delta $  \hspace{2ex} and \hspace{2ex}    $\displaystyle 1 + \frac{\delta}{2}< 2-\frac{m\gamma}{N_m} < 1+\delta. $

\vspace{2ex}

\noindent Since $0 < 2r < 1,$ there exists a constant $C_2(r) > 0$ depending only on $r,$ such that $C(r, \, \frac{m\gamma}{N_m}) \leq C_2(r),$ for all $m\in \mathbb{N}.$ Since $r_j \leq 2r < 1$, we have:

\vspace{2ex}

\hspace{6ex}  $r_j^{2(2-\frac{m\gamma}{N_m})} < r_j^2,$ \, for \hspace{1ex} $2(2-\frac{m\gamma}{N_m}) > 2.$

\vspace{2ex}

\noindent Thus, estimate $(e')$ (inferred above from $(e)$ of Lemma \ref{Lem:atomisation}) implies:

$$ \int_{D}|f_m(z)|^2 \,  e^{-2m\varphi_0(z)} \, d\lambda (z) \leq C(r), \hspace{2ex} \forall m\gg 0,$$

\noindent where $C(r)=  C_1(r) \, C_2(r) > 0$ is a constant depending only on the radius $r$ of the disc $D$ on which we are working. This yields conclusion $(iii)$ of Theorem \ref{The:potentiel} and completes its proof.   \hfill $\Box$

\section{Preliminaries to Theorem \ref{The:globalregularisation}}\label{subsection:local}

In this section, we recall some classical material, mainly from [Dem92], that will be used throughout the paper.

\subsection{The Demailly regularisation of psh functions}\label{subsection:Demailly_psh}

  Let $\varphi$ be a psh function on a bounded pseudoconvex open set $\Omega \subset \C^n$. A well-known result of Demailly (cf. [Dem92, Proposition 3.1]) asserts that $\varphi$ can be approximated pointwise and in $L^1_{loc}(\Omega)$ topology by psh functions $\varphi_m$ with analytic singularities (see definition (\ref{equation:analsing}) in the introduction), defined as 

\begin{equation}\label{eqn:Demaillyapprox}\varphi_m= \frac{1}{2m}\, \log \sum\limits_{j=0}^{+\infty} |\sigma_{m, \, j}|^2,\end{equation}

\noindent where $(\sigma_{m, \, j})_{j\in \mathbb{N}}$ is an arbitrary orthonormal basis of the Hilbert space ${\mathcal H}_{\Omega}(m\varphi)$ of holomorphic functions $f$ on $\Omega$ such that $|f|^2\, e^{-2m\varphi}$ is integrable on $\Omega$. They even satisfy the estimates:

\begin{equation}\label{eqn:approximation}\varphi(z) - \frac{C_1}{m} \leq \varphi_m(z) \leq \sup_{|\zeta - z|<r}\varphi(\zeta) + \frac{1}{m} \, \log \, \frac{C_2}{r^n},\end{equation}

\noindent for every $z\in \Omega$ and every $r< d(z, \, \partial \Omega)$. In particular, the sequence $dd^c\varphi_m$ converges to $dd^c\varphi$ in the weak topology of currents, and the corresponding Lelong numbers satisfy:

\begin{equation}\label{eqn:Lelongnumbers}\nu(\varphi, \, x)-\frac{n}{m} \leq \nu(\varphi_m, \, x) \leq \nu(\varphi, \, x), \hspace{2ex} x\in \Omega.\end{equation}

\noindent For analytic singularities, the Lelong number $\nu(\varphi_m, \, x)$ at an arbitrary point $x$ equals $\frac{1}{m}\, \min\limits_{j\geq 0} \mbox{ord}_x \sigma_{m, \, j}$, where $\mbox{ord}_x$ is the vanishing order at $x$. The sequence $(\varphi_m)_{m\in \N}$ defined in $(\ref{eqn:Demaillyapprox})$ has come to be referred to as the {\bf Demailly approximation} of $\varphi$.

\subsection{An example}\label{subsection:example} In order to give a taster of the difficulties involved in the proof of Theorem \ref{The:globalregularisation}, consider for a moment the case when $\varphi$ has a zero Lelong number at every point $x\in \Omega$. It may seem tempting to alter the Demailly regularisation of $\varphi$ by defining the sequence of smooth psh functions $(\psi_m)_{m\in \mathbb{N}}$ on $\Omega$ as \begin{equation}\label{eqn:masses} \psi_m:=\frac{1}{2m}\, \log\bigg(\sum\limits_{j=0}^{+\infty} |\sigma_{m, \, j}|^2 + \sum\limits_{j=0}^{+\infty} \bigg| \frac{\partial\sigma_{m, \, j}}{\partial z_1}\bigg|^2 + \dots + \sum\limits_{j=0}^{+\infty} \bigg| \frac{\partial\sigma_{m, \, j}}{\partial z_n}\bigg|^2\bigg),\end{equation} \noindent where $\frac{\partial}{\partial z_1}, \dots ,\, \frac{\partial}{\partial z_n}$ are the first order partial derivarives with respect to the standard coordinates $z=(z_1, \dots , \, z_n)$ on $\C^n$. Then $dd^c\psi_m$ converges to $dd^c\varphi$ in the weak topology of currents as $m\rightarrow +\infty$. However, all we can get from an application of Theorem \ref{The:potentiel} and of the Ohsawa-Takegoshi $L^2$ extension theorem is the following upper estimate of the Monge-Amp\`ere masses on every relatively compact open subset $B\Subset \Omega$\!\!:

\hspace{10ex} $\displaystyle \int\limits_{B}(dd^c\psi_m)^k \wedge \beta^{n-k} \leq C_m\, (\log m)^k$, \hspace{2ex} $k=1, \dots , \, n$,  \\

\noindent where $\beta$ is the standard K\"ahler form on $\C^n$ and $C_m>0$ is a constant depending on $m$ in a way that cannot be controlled at this stage\footnote{An earlier version of this work claimed that the constant $C_m$ was independent of $m$. That mistake did not affect in any way either the statements or the proofs of Theorems \ref{The:potentiel}, \ref{The:globalregularisation} and \ref{The:bigness} and is now corrected. The functions $\psi_m$ were correctly estimated below with uniform constants at certain points $a_{m,\,j}$, but that estimate did not extend to all the points.}. 

\vspace{2ex}

 It is worth stressing that the sequence of Monge-Amp\`ere masses may not be bounded above. To see this, suppose $\varphi$ is $C^{\infty}$ in the complement of an analytic set $V\subset \Omega$. If $(\varphi_m)_{m\in \N}$ is the Demailly approximation of $\varphi$, it is shown in [DPS01, p. 701-702] that the sequence $(\varphi_{2^m}+ 2^{-m})_{m\in \N}$ is decreasing using an effective version of the subadditivity property of multiplier ideal sheaves. The same proof shows that the corresponding sequence $(\psi_{2^m}+ 2^{-m})_{m\in \N}$ in the regularisation defined in (\ref{eqn:masses}) is also decreasing. Then the $C^{\infty}$ $(1, \, 1)$-forms $(dd^c\psi_{2^m})^k$ are well-defined on $\Omega \setminus V$ and converge in the weak topology of currents to $(dd^c\varphi)^k$ for every $k=1, \dots , \, n$ (see [Dem97, chapter III, Theorem 3.7]). If $K\Subset B\setminus V$ and $0\leq \chi \leq 1$ is a $C^{\infty}$ function with compact support in $B\setminus V$ such that $\chi \equiv 1$ on $K$, then:\\

$\displaystyle \int_{B\setminus V} \chi (dd^c\psi_{2^m})^k \wedge \beta^{n-k} \leq \int_{B\setminus V} (dd^c\psi_{2^m})^k \wedge \beta^{n-k}, \hspace{2ex} m\in \N$, \\

\noindent and taking $\liminf\limits_{m\rightarrow +\infty}$, the weak convergence implies: \\

$\displaystyle \int_{B\setminus V} \chi (dd^c\varphi)^k \wedge \beta^{n-k} \leq \liminf_{m\rightarrow +\infty}\int_{B\setminus V} (dd^c\psi_{2^m})^k \wedge \beta^{n-k}, \hspace{2ex} k=1, \dots , \, n$. \\

\noindent Clearly $\int_K (dd^c\varphi)^k \wedge \beta^{n-k} \leq \int_{B\setminus V} \chi (dd^c\varphi)^k \wedge \beta^{n-k}$, and letting $K\Subset B\setminus V$ increase, we get: \\

$\displaystyle \int_{B\setminus V} (dd^c\varphi)^k \wedge \beta^{n-k} \leq \liminf_{m\rightarrow +\infty}\int_{B\setminus V} (dd^c\psi_{2^m})^k \wedge \beta^{n-k}, \hspace{2ex} k=1, \dots , \, n$. \\

\noindent Now, there are examples of psh functions $\varphi$ for which the Monge-Amp\`ere mass on the left-hand side above is infinite for $k=n$ (see Kiselman's example in [Kis84, p.141-143] of a $\varphi$ with zero Lelong numbers, or the Shiffman-Taylor example in [Siu 75, p.451-453]). Thus the last inequality shows that for such functions the sequence of Monge-Amp\`ere masses associated with the above regularisation is unbounded.

\subsection{The patching procedure}\label{subsection:global}

 In this subsection, we recall how local regularisations can be patched together  by slightly modifying the classical patching procedure introduced in [Dem92] to produce a global regularisation of a (quasi-)positive closed $(1,\,1)$-current. For the sake of simplicity, we assume $X$ to be compact. The method actually works for any manifold $X$ which can be covered by finitely many coordinate patches on which the local regularisation procedure can be applied. We will explain the procedure in the case of zero Lelong numbers for the sake of simplicity. 

\vspace{2ex}

 Let $T\geq \gamma$ be a $d$-closed current of bidegree $(1, \, 1)$ on a compact complex manifold $X$, where $\gamma$ is a continuous $(1, \, 1)$-form such that $d\gamma=0$. Assume $T$ has a zero Lelong number at every point in $X$. Then, there exist $C^1$ $(1, \, 1)$-forms $T_m$ in the same $\partial\bar{\partial}$-cohomology class as $T$ which converge to $T$ in the weak topology of currents and satisfy\!\!: \\

\noindent $(a)$ \, $\displaystyle T_m \geq \gamma - \frac{C}{m}\, \omega$; \\

\noindent $(b)$ \,  $\displaystyle \int\limits_{X}(T_m -\gamma + \frac{C}{m}\, \omega)^q \wedge \omega^{n-q} \leq C_m\, (\log m)^q$, \hspace{2ex} $q=1, \dots , \, n=\mbox{dim}_{\C}X$, \\

\noindent for a fixed Hermitian metric $\omega$ on $X$ and some constant $C_m>0$ depending on $m$ in a way that cannot be controlled at this stage.

\vspace{2ex}

 As the patching procedure is essentially well-known (see, for instance, [Dem92, section 3]), we will only point out the new aspects arising when Monge-Amp\`ere masses are being estimated.

 We always have $T=\alpha + dd^c\varphi \geq \gamma$ globally on $X$, with some $C^{\infty}$ $(1, \, 1)$-form $\alpha$ and some almost psh function $\varphi$ on $X$. The set-up is the one described in the introduction. After possibly replacing $T$ with $T-\alpha$ and $\gamma$ with $\gamma -\alpha$, we can assume
$T= dd^c\varphi \geq \gamma$. Let us fix $\delta >0$, and four finite coverings of $X$ by concentric coordinate balls $(B^{(3)}_j)_j,$ $(B'_j)_j,$ $(B''_j)_j$ and $(B_j)_j$ of radii $\frac{\delta}{2},$ $\delta,$ $\frac{3}{2}\delta$, and respectively $2\delta.$ Since $d\gamma =0$, $\gamma$ is locally exact and we can assume that, for every $j$, $\gamma = dd^ch_j$ on $B_j$ for some $C^1$ function $h_j$. The function: \\

\hspace{6ex} $\psi_j: = \varphi - h_j,$  \\

\noindent is psh on $B_j$ for every $j$. As, by assumption, $\psi_j$ has a zero Lelong number at every point in $B_j$, construction (\ref{eqn:masses}) can be applied to each $\psi_j$ on $B_j$ to get smooth approximations: \\

\hspace{6ex} $\displaystyle \psi_{j, \, m}: = \frac{1}{2m}\, \log\bigg(\sum\limits_{l=0}^{+\infty} |\sigma_{j, \, m, \, \, l}|^2  + \sum\limits_{r=1}^n \sum\limits_{l=0}^{+\infty} \bigg|\frac{\partial\sigma_{j,\, m,\, l}}{\partial z_r}\bigg|^2\bigg),$ \\

\noindent with an arbitrary orthonormal basis $(\sigma_{j, \, m, \, \, l})_{l\in \N}$ of the Hilbert space ${\mathcal H}_{B_j}(m\psi_j)$ (see notation in the previous section). Then $\varphi_{j, \, m}:= \psi_{j, \, m} + h_j$ converges pointwise and in  $L^1_{loc}$ topology to $\varphi$ as $m\rightarrow +\infty$ on $B_j$, and these local approximations can be glued together into a global approximation of $\varphi$ defined as:  \\

\hspace{6ex}$\displaystyle \varphi_m(z): = \sup\limits_{B''_j \ni z} \bigg(\varphi_{j, \, m}(z) + \frac{C_1(\delta)}{m}\, (\delta^2 - |z^j|^2)\bigg),$  \\

\noindent with a constant $C_1(\delta)>0$ depending only on $\delta$ which will be specified below, and a local holomorphic coordinate system $z^j$ centred at the centre of $B_j$. The currents $T_m:=dd^c\varphi_m$ satisfy the conclusions of Theorem \ref{The:zerolreg} if the following patching condition holds:

\begin{equation}\label{eqn:patching}\displaystyle \varphi_{j, \, m}(z) + \frac{C_1(\delta)}{m}\, (\delta^2 - |z^j|^2) \leq \varphi_{k, \, m}(z) + \frac{C_1(\delta)}{m}\, (\delta^2 - |z^k|^2),\end{equation}

\noindent for $z\in (\overline{B}''_j \setminus B'_j)\cap B^{(3)}_k.$ One can then prove the existence of a constant $C_1(\delta)>0$ satisfying this patching condition by means of H\"ormander's $L^2$ estimates ([Hor65]). One need only estimate the difference $\psi_{j, \, m} -  \psi_{k, \, m}$ on $B''_j \cap B''_k$ and show that $\varphi_{j, \, m} -  \varphi_{k, \, m}$ is uniformly bounded above on $B''_j \cap B''_k$ by $O(\frac{1}{m})$ as $m\rightarrow +\infty$. Now, for every fixed $z\in B_j$, the norms of the linear maps $f\mapsto f(z)$ and $f\mapsto \frac{\partial f}{\partial z_r}(z)$, $r=1, \dots ,\, n,$ defined on the Hilbert space ${\mathcal H}_{B_j}(m\psi_j)$, can be expressed in terms of an orthonormal basis, and we get: \\

$\displaystyle \frac{1}{2m}\, \log \sup\limits_{f\in\overline{B}_{j, \, m}} \bigg(|f(z)|^2 + \sum\limits_{r=1}^n \bigg |\frac{\partial f}{\partial z_r}(z)\bigg |^2 \bigg ) \leq \psi_{j,\, m}(z) $ \\

 \hspace{20ex}  $\displaystyle \leq   \frac{1}{2m}\, \log\bigg( (n+1) \, \sup\limits_{f\in\overline{B}_{j, \, m}}\bigg(|f(z)|^2 + \sum\limits_{r=1}^n \bigg |\frac{\partial f}{\partial z_r}(z)\bigg |^2 \bigg )\bigg)$, \\

\noindent where $\overline{B}_{j, \, m}$ is the unit ball of ${\mathcal H}_{B_j}(m\psi_j)$. We also have the analogous relations for $\psi_{k,\, m}$ on $ B_k$. This means that to compare $\psi_{j,\, m}$ and $\psi_{k,\, m}$ at a fixed point $x_0\in B''_j \cap B''_k$, it is enough to show that for every holomorphic function $f_j$ on $B_j$ such that $\int_{B_j}|f_j|^2 \, e^{-2m\psi_j} =1$, there exists a holomorphic function $f_k$ on $B_k$ having an $L^2$-norm under control and satisfying: \\

$f_k(x_0)=f_j(x_0)$,  \hspace{2ex} and \hspace{2ex}  $\displaystyle \frac{\partial f_k}{\partial z_r}(x_0) = \frac{\partial f_j}{\partial z_r}(x_0)$, \hspace{2ex} for $r=1, \dots , \, n$. \\

\noindent This is done using H\"ormander's $L^2$ estimates ([Hor65]). Let $\theta$ be a cut-off function supported in a neighbourhood of $x_0$ such that $\theta \equiv 1$ near $x_0$, and solve the equation 

\begin{equation}\label{eqn:horeqn}\bar{\partial} g = \bar{\partial} (\theta f_j)\end{equation}

\noindent on $B_k$ with a weight containing the term $2(n+1)\, \log|z-x_0|$ which forces the solution $g$ to vanish to order at least $2$ at $x_0$. Specifically, if $h_{jk}$ is a holomorphic function on $B_j \cup B_k$ such that $h_j - h_k =\mbox{Re}\, h_{jk}$ on $B_j \cap B_k$, we can find a solution $g$ to the above equation on $B_k$ satisfying H\"ormander's $L^2$ estimates with the strictly psh weight: 

\begin{equation}\label{eqn:weight}2m (\psi_k - \mbox{Re}\, h_{jk}) + 2(n+1)\, \log|z-x_0| + |z-x_0|^2.\end{equation}

\noindent Now set $f_k:=\theta f_j - g$ which is easily seen to satisfy the requirements. The precise estimate of the solution $g$ gives the uniform upper estimate of $\varphi_{j, \, m}- \varphi_{k, \, m}$ on $B''_j \cap B''_k$ by $O(\frac{1}{m})$ which implies the existence of a constant $C_1(\delta)>0$ satisfying the patching condition (\ref{eqn:patching}). The details are left to the reader.

 The loss of positivity incurred in $T_m$ with respect to the original $T$ can be seen to be at most $\frac{C}{m}$ as in [Pop04] thanks to the form $\gamma$ being closed. This proves $(a)$. That the approximating currents $T_m:=dd^c\varphi_m$ constructed through this patching procedure satisfy the condition $(b)$ on Monge-Amp\`ere masses follows from the local construction (\ref{eqn:masses}) of the previous subsection.

\section{Modified regularisation of currents}\label{subsection:modreg}

 In this section, we start the main construction of the paper.

 Let $\Omega\subset \C^n$ be a bounded pseudoconvex open set, and let $\varphi$ be a strictly psh function on $\Omega$ such that $i\partial\bar{\partial}\varphi \geq C_0\, \beta$ for some constant $C_0>0$ and the standard K\"ahler form $\beta$ on $\C^n$. Our goal is to construct regularising psh functions with analytic singularities approximating $\varphi$ for which the Monge-Amp\`ere masses can be controlled. The overall idea is to modify Demailly's regularisation (\ref{eqn:Demaillyapprox}) by adding derivatives of the functions $\sigma_{m, \, j}, j\in \N,$ which form an orthonormal basis of the Hilbert space:

\begin{equation}\label{eqn:Hilbdef}\displaystyle{\mathcal H}_{\Omega}(m\varphi):=\bigg\{f\in{\mathcal O}(\Omega)\,\, ; \,\, \int\limits_{\Omega}|f|^2\, e^{-2m\varphi}\, dV_n < +\infty\bigg\}, \hspace{2ex} dV_n:=\frac{\beta^n}{n!}.\end{equation}

\noindent We shall repeatedly denote $D^{\alpha}$ the derivation operator with respect to any multi-index $\alpha=(\alpha_1, \dots , \alpha_n)\in\N^n$ of length $|\alpha|=\alpha_1+ \dots +\alpha_n$, and $[\,\,\,]$ the integer part. Unlike section \ref{subsection:local} where it was enough to derive to order one as $\varphi$ was assumed to have zero Lelong numbers everywhere, we need to derive more in the general case. Yet we cannot afford to derive too much and have to accept singularities in the regularising functions as shown below.

\begin{Lem}\label{Lem:derivorder} For every $\delta>0$, the psh functions with analytic singularities defined as: 

\hspace{6ex}$\displaystyle \varphi_m^{\delta} = \frac{1}{2m} \log \sum\limits_{j=0}^{+\infty}\sum\limits_{|\alpha|=0}^{[\delta m]} \bigg| \frac{D^{\alpha} \sigma_{m, \, j}}{\alpha !}  \bigg|^2$, \hspace{2ex} $m\in\N,$ \\

\noindent satisfy the estimates below for constants $C_1, C_3>0$ independent of $m$ and $\varphi$: \\

$\displaystyle\varphi(z) - \frac{C_1}{m} \leq \varphi_m^{\delta}(z) \leq \sup\limits_{|\zeta - z|< 2r} \varphi(\zeta) - \bigg(\frac{[\delta m]}{m}+\frac{n}{m}\bigg)\, \log r + \frac{1}{m} \log C_3,$ \\

\noindent at every point $z\in \Omega$ and for every $0\leq r < \min\{\frac{1}{2} d(z, \, \partial \Omega), \, 1\}$. In particular, if we choose $\delta_m:=C\, \varepsilon_m$ with $\varepsilon_m\downarrow 0$ and $C>0$ a constant independent of $m$, $\varphi_m^{\delta_m}$ converges pointwise and in $L^1_{loc}$ topology to $\varphi$ when $m\rightarrow +\infty$.

\end{Lem}

\noindent {\it Proof.} The lower bound follows from the lower estimate in (\ref{eqn:approximation}) since $\varphi_m^{\delta} \geq \varphi_m$. To get the upper bound, we apply Parseval's formula to each function $\sigma_{m, \, j}$ on the sphere $S(z, \, r)$ and then sum over $j$ to get: \\

$\begin{array}{lll}

\frac{C}{r^{2n-1}} \, \displaystyle \int_{S(z,r)}\overset{+\infty}{\underset{j=0}\sum} |\sigma_{m,j}(\zeta)|^2 \, d\sigma(\zeta) & = & \displaystyle \sum\limits_{j=0}^{+\infty}\sum\limits_{\alpha\in \mathbb{N}^n} \bigg |\frac{D^{\alpha}\sigma_{m,j}}{\alpha!}(z)\bigg |^2 \, r^{2|\alpha|} \\

  & \geq & r^{2[\delta m]} \displaystyle \sum\limits_{j=0}^{+\infty}\sum\limits_{|\alpha|=0}^{[\delta m]} \bigg |\frac{D^{\alpha}\sigma_{m,j}}{\alpha!}(z) \bigg|^2.

\end{array}$ \\

\noindent It is now enough to take $\frac{1}{2m}\, \log$ on both sides and to use the upper estimate in (\ref{eqn:approximation}) to conclude. \hfill $\Box$

\vspace{2ex}

 In other words, we still have an approximation of $\varphi$ if we derive the $\sigma_{m, \, j}$'s up to order $\leq o(m)$ (e.g. $[C\varepsilon_m m]$ for some $\varepsilon_m \downarrow 0$, even when $[C\varepsilon_m m]\rightarrow +\infty$). But we would lose control on the upper bound of $\varphi_m^{\delta}$ if we derived up to removing all the common zeroes of the $\sigma_{m, \, j}$'s (e.g. up to order $[Am]$ with $A:=\sup\limits_{x\in \Omega}\nu(\varphi, \, x)$ assumed to be finite after possibly shrinking $\Omega$). \\

 Besides taking derivatives, we will further modify Demailly's regularising functions (\ref{eqn:Demaillyapprox}). A crucial role will be played by multiplier ideal sheaves ${\mathcal I}(m\varphi)\subset{\mathcal O}_{\Omega}$ associated with psh functions. They are defined, for every $m\in\N$, as (cf. [Nad90], [Dem93]):

\begin{equation}\label{eqn:misdef} {\mathcal I}(m\varphi)_x:=\{f\in {\mathcal O}_{\Omega,\, x}, \,\, |f|^2\, e^{-2m\varphi} \,\, \mbox{is Lebesgue-integrable near}\,\, x\},\end{equation}

\noindent at every $x\in\Omega$. It is a well-known result of Nadel ([Nad90], see also [Dem93]) that, for every $m\in \N$, the multiplier ideal sheaf ${\mathcal I}(m\varphi)$ is coherent and is generated as an ${\mathcal O}_{\Omega}$-module by an arbitrary orthonormal basis $(\sigma_{m, \, j})_{j\in \N}$ of the Hilbert space ${\mathcal H}_{\Omega}(m\varphi)$. By coherence, the restriction of ${\mathcal I}(m\varphi)$ to every compact subset has only finitely many generators $(\sigma_{m, \, j})_{1\leq j\leq N_m}$. When $m\rightarrow +\infty$, this local finite generation property was made effective in [Pop06, Theorem 1.1] in the following way. It was shown there that for any given relatively compact open subset $B\Subset \Omega$, there exist a subset $B_0\Subset B$ and $m_0=m_0(C_0)\in \N$ such that for every $m\geq m_0$ one can find an orthonormal basis $(\sigma_{m, \, j})_{j\in \N}$ of ${\mathcal H}_{\Omega}(m\varphi)$ and finitely many elements $\sigma_{m, \, 1}, \dots , \sigma_{m, \, N_m}$ in it with the following property. Every local section $g\in {\mathcal H}_B(m\varphi)$ admits a decomposition:

\begin{equation}\label{eqn:redecomp}g(z)=\sum\limits_{j=1}^{N_m}h_{m, \, j}(z)\, \sigma_{m, \, j}(z),  \hspace{3ex} z\in B_0,\end{equation}

\noindent with some holomorphic functions $h_{m, \, j}$ on $B_0$ satisfying: 

\begin{equation}\label{eqn:reestimate}\sup\limits_{B_0}\sum\limits_{j=1}^{N_m}|h_{m, \, j}|^2 \leq C\, N_m\, \int\limits_B |g|^2\, e^{-2m\varphi} < +\infty,\end{equation}

\noindent for a constant $C>0$ depending only on $n,\, r,$ and the diameter of $\Omega$. Furthermore, in the case of a $\varphi$ with analytic singularities, the growth of the number $N_m$ of local generators needed is at most polynomial of degree $n$ as $m\rightarrow +\infty$, namely

\begin{equation}\label{eqn:polgrowth}\lim\limits_{m\rightarrow +\infty}\frac{n!}{m^n}\, N_m = \frac{2^n}{\pi^n} \, \int\limits_B (i\partial\bar{\partial}\varphi)^n_{ac} < +\infty,\end{equation}

\noindent where $(i\partial\bar{\partial}\varphi)_{ac}$ denotes the absolutely continuous part of $i\partial\bar{\partial}\varphi$ in the Lebesgue decomposition of its measure coefficients into an absolutely continuous and a singular part with respect to the Lebesgue measure. An immediate corollary of formula (\ref{eqn:redecomp}) and estimate (\ref{eqn:reestimate}) is the following comparison relation:

\begin{equation}\label{eqn:recomparison}\sum\limits_{j=0}^{+\infty}|\sigma_{m, \, j}|^2 \leq C\, N_m \, \sum\limits_{j=1}^{N_m}|\sigma_{m, \, j}|^2 \hspace{3ex} \mbox{on} \,\, B_0.\end{equation} 
   
\noindent Building on these results from [Pop06], we can define a new regularisation of $\varphi$ using only finitely many generators of ${\mathcal I}(m\varphi)_{|B_0}$ to modify the definition in the previous Lemma \ref{Lem:derivorder}.

\begin{Lem}\label{Lem:newreg} The psh functions with analytic singularities: \\

\hspace{6ex}$\displaystyle \psi_m = \frac{1}{2m} \log \sum\limits_{j=1}^{N_m}\sum\limits_{|\alpha|=0}^{[C\varepsilon_m m]} \bigg| \frac{D^{\alpha} \sigma_{m, \, j}}{\alpha !}  \bigg|^2,$ \hspace{2ex} $m\in \N,$ \\

\noindent satisfy the estimates below for constants $C_1, C_3>0$ independent of $m$ and $\varphi$: \\

\noindent $\displaystyle\varphi(z) - \frac{\log N_m +C_1}{2m} \leq \psi_m(z) \leq \sup\limits_{|\zeta - z|< 2r} \varphi(\zeta) - \frac{[C\varepsilon_m m]+n}{m}\, \log r + \frac{1}{m} \log C_3,$ \\

\noindent at every point $z\in B_0$ and for every $0\leq r < \min\{\frac{1}{2}\, \mbox{dist}\,(z, \, \partial \Omega), \, 1\}$. In particular, if $\varepsilon_m\downarrow 0$ and $C>0$ is a constant independent of $m$, $\psi_m$ converges pointwise on $B_0$ and in $L^1_{loc}(B_0)$ topology to $\varphi$ when $m\rightarrow +\infty$.

\end{Lem}

\noindent {\it Proof.} It follows immediately from the above considerations. \hfill $\Box$ \\

 This regularisation is not yet satisfactory from the point of view of the Monge-Amp\`ere masses as $\psi_m$ still has singularities and the Chern-Levine-Nirenberg inequalities cannot be applied as such (see section \ref{subsection:local} where they were applied to smooth psh functions). Building on ideas in this section, we will construct a new regularisation of $\varphi$ in subsequent sections for which we can control the Monge-Amp\`ere masses.

\section{Additivity defect of multiplier ideal sheaves}\label{subsection:addefect}

 Throughout this section we shall suppose that $\varphi$ is a psh function with analytic singularities on $\Omega\Subset \C^n$ of the form:

\begin{equation}\label{eqn:assumption}\displaystyle \varphi =\frac{c}{2}\, \log(|g_1|^2 + \dots + |g_N|^2) + v,\end{equation}

\noindent for some (possibly infinitely many) holomorphic functions $g_1, \dots , g_N\in {\mathcal O}(\Omega)$, some constant $c>0$, and some $C^{\infty}$ function $v$ on $\Omega$. We also assume that $i\partial\bar{\partial}\varphi\geq C_0\, \beta$ for some constant $C_0>0$ and the standard K\"ahler form $\beta$ on $\C^n$. The results in this section will be subsequently applied to $\varphi_p$ in place of $\varphi$, where $(\varphi_p)_{p\in \N}$ are the Demailly regularisations (\ref{eqn:Demaillyapprox}) of an arbitrary psh function $\varphi$. We work with general functions of the form (\ref{eqn:assumption}) in this section for the sake of generality.

Fix now a relatively compact pseudoconvex open subset $B\Subset \Omega$. Our aim is to find better regularisations of $\varphi$ for which we can control the Monge-Amp\`ere masses. In so doing, we will still use the ideal sheaves ${\mathcal I}(m\varphi)$, but the main idea underlying the argument is to take $m=m_0q$ with $m\gg m_0$ (though $m_0\rightarrow +\infty$) and to reduce the study of ${\mathcal I}(m\varphi)$ to the study of ${\mathcal I}(m_0\varphi)$. We shall then get a grip on ${\mathcal I}(m\varphi)$ as being ``almost'' equal to ${\mathcal I}(m_0\varphi)^q$ up to an error that we are able to estimate. This is made precise in the following.

\begin{Prop}\label{Prop:inclusions} For any $\varepsilon >0$, any $m_0\geq \frac{n+2}{c\,\varepsilon}$, any $q \in \N$, and any $B\Subset \Omega$, the following inclusions of multiplier ideal sheaves hold:

\begin{equation}\label{eqn:inclusions}{\mathcal I}(m_0(1+\varepsilon)\varphi)_{|B}^q \subset {\mathcal I}(m_0q\varphi)_{|B} \subset {\mathcal I}(m_0\varphi)_{|B}^q.\end{equation}

\end{Prop}

\noindent The right-hand inclusion actually holds on $\Omega$ for every $m_0$ and is the subadditivity property of multiplier ideal sheaves proved by Demailly, Ein and Lazarsfeld in [DEL00]. It relies on the Ohsawa-Takegoshi $L^2$ extension theorem ([OT87]). The left-hand inclusion was proved in [Pop06] using Skoda's $L^2$ division theorem ([Sko72b]). It can be seen as measuring the extent to which multiplier ideal sheaves fail to have an additive growth. These inclusions can be given effective versions, with estimates, that we are now undertaking to make explicit.

 Consistent with the notation in (\ref{eqn:Hilbdef}), consider the Hilbert spaces: \\

\hspace{6ex}$\displaystyle {\mathcal H}_{\Omega}(m_0(1+\varepsilon)\, \varphi)$, $\displaystyle {\mathcal H}_{\Omega}(m_0q\, \varphi)$, and $\displaystyle {\mathcal H}_{\Omega}(m_0\, \varphi)$, \\

\noindent and respective orthonormal bases: \\

\hspace{6ex} $\displaystyle (\sigma_{m_0(1+\varepsilon), \, j})_{j\in \N}$, $\displaystyle (\sigma_{m_0q, \, j})_{j\in \N}$, and $\displaystyle (\sigma_{m_0, \, j})_{j\in \N}$. \\

\noindent The multiplier ideal sheaves in (\ref{eqn:inclusions}) are generated as ${\mathcal O}_{\Omega}$-modules by the above Hilbert space orthonormal bases respectively. Since they are coherent, their restrictions to the relatively compact subset $B$ are finitely generated. Possibly after reordering, we can therefore assume that they are generated as follows: \\

$\displaystyle {\mathcal I}(m_0(1+\varepsilon)\, \varphi)_{|B}=(\sigma_{m_0(1+\varepsilon), \, j})_{1\leq j\leq N_{m_0(1+\varepsilon)}}$, \\

 $\displaystyle {\mathcal I}(m_0q \, \varphi)_{|B}=(\sigma_{m_0q, \, j})_{1\leq j\leq N_{m_0q}}$, \hspace{2ex} $\displaystyle {\mathcal I}(m_0 \, \varphi)_{|B}= (\sigma_{m_0, \, j})_{1\leq j\leq N_{m_0}}$. \\ 

\noindent The results of [Pop06], summarised in the previous section as (\ref{eqn:redecomp}) -- (\ref{eqn:recomparison}), will be made an essential use of in all that follows. We start by stating the following effective versions of the left-hand (cf. $(a)$) and right-hand (cf. $(b)$) inclusions in (\ref{eqn:inclusions}) when derivatives are also taken into account.

\begin{Prop}\label{Prop:effectiveincl} 

\noindent $(a)$ Let $\varphi$ be a psh function on $\Omega\Subset \C^n$ of the form (\ref{eqn:assumption}). Given $B\Subset\Omega$, there exist an open subset $B_0\Subset B$ and an orthonormal basis $(\sigma_{m_0q, \, j})_{j\in \N}$ of ${\mathcal H}_{\Omega}(m_0q\, \varphi)$ such that $\sigma_{m_0q, \, j},$ $j=1, \dots, N_{m_0q},$ generate the ${\mathcal O}_{\Omega}$-module ${\mathcal I}(m_0q\, \varphi)$ on $B_0$ and the following estimate holds at every $z\in B_0$ for every multi-index $\alpha=(\alpha_1, \dots , \alpha_n)\in\N^n$ and every orthonormal basis $(\sigma_{m_0(1+\varepsilon), \, j})_{j\in \N}$ of ${\mathcal H}_{\Omega}(m_0(1+\varepsilon)\varphi):$ \\

\noindent \begin{equation}\label{eqn:superad}\displaystyle \sum\limits_{j_1, \dots , j_q=0}^{+\infty} | D^{\alpha} (\sigma_{m_0(1+\varepsilon),\,j_1} \dots \sigma_{m_0(1+\varepsilon),\,j_q})|^2 \leq  \displaystyle C_{\alpha}\, N_{m_0q}^2 \, C_{m_0}^q \, \sum\limits_{j=1}^{N_{m_0q}}\sum\limits_{\beta\leq \alpha}|D^{\beta}\sigma_{m_0q, \, j}|^2,\end{equation} 

\vspace{2ex}

\noindent for any $q \in \N$, any $0 < \varepsilon << 1$, and any $m_0\geq \frac{n+2}{c\,\varepsilon}$. The constants: \\

  $0<C_{\alpha} \leq O\bigg(|\alpha|^2\, \sum\limits_{\beta\leq \alpha}{\alpha \choose \beta}^2\bigg)$, \hspace{2ex} $C_{m_0}:= C_n(m_0\,c\,(1+\varepsilon)-n)\, (\sup\limits_Be^{\varphi})^{2m_0\varepsilon}$ \\

\noindent are such that $C_{\alpha}$ depends only on $\alpha$, $n$, $B_0,$ $B$, and $C_n>0$ depends only on $n,\, B,\, \Omega$; \\

\noindent $(b)$ For any orthonormal bases of ${\mathcal H}_{\Omega}(m_0q\varphi)$ and ${\mathcal H}_{\Omega}(m_0\varphi)$, we have:

\begin{equation}\label{eqn:subad}\displaystyle \sum\limits_{j=0}^{+\infty}|D^{\alpha}\sigma_{m_0q, \, j}|^2 \leq C_n^{q-1}\, \sum\limits_{j_1, \dots , j_q=0}^{+\infty}|D^{\alpha}(\sigma_{m_0, \, j_1} \dots \sigma_{m_0,\,j_q})|^2 \hspace{3ex} {\rm on} \hspace{1ex} \Omega,\end{equation} 

\noindent for every multi-index $\alpha=(\alpha_1, \dots , \alpha_n)\in\N^n$, with a constant $C_n>0$ depending only on $n$ and $\Omega$.

\end{Prop}

\noindent {\it Proof.} $(b)$ \, Estimate (\ref{eqn:subad}) follows from the effective version of the subadditivity property of multiplier ideal sheaves established in [DEL00] (see also [DPS01, proof of Theorem 2.2.1.]). Indeed, for every $f\in {\mathcal O}(\Omega)$ satisfying \\

\hspace{6ex}$\displaystyle \int\limits_{\Omega}|f|^2\, e^{-2q\, (m_0\, \varphi)}\, dV_n=1,$ \\

\noindent Step 3 in the proof of Theorem 2.2.1. in [DPS01], when applied in this context, shows by means of the Ohsawa-Takegoshi $L^2$ extension theorem applied $q-1$ times that, at every point $z\in \Omega$, there is a decomposition: 

\begin{equation}\label{eqn:subaddecomp}f(z)= \sum\limits_{j_1, \dots , j_q=0}^{+\infty}c_{j_1, \dots , j_q} \, \sigma_{m_0, \, j_1}(z) \dots \sigma_{m_0, \, j_q}(z), \hspace{2ex} \sum\limits_{j-1, \dots, j_q=0}^{+\infty}|c_{j_1, \dots , j_q}|^2 \leq C_n^{q-1},\end{equation}

\noindent with scalar coefficients $c_{j_1, \dots , j_q}\in \C$ satisfying the above estimate for some constant $C_n>0$ depending only on $n$ and $\Omega$. Applying $D^{\alpha}$ we get: \\

\hspace{6ex}$\displaystyle D^{\alpha}f(z) = \sum\limits_{j_1, \dots , j_q=0}^{+\infty}c_{j_1, \dots , j_q}\, D^{\alpha}(\sigma_{m_0,\,j_1} \dots \sigma_{m_0,\,j_q})(z),$ \hspace{3ex} $z\in \Omega,$ \\

\noindent The Cauchy-Schwarz inequality and estimate (\ref{eqn:subaddecomp}) of the scalar coefficients $c_{j_1, \dots , j_q}$ give: \\

\hspace{6ex} $\displaystyle |D^{\alpha}f(z)|^2 \leq C_n^{q-1}\, \sum\limits_{j_1, \dots , j_q=0}^{+\infty}|D^{\alpha}(\sigma_{m_0, \, j_1}\dots \sigma_{m_0, \, j_q})(z)|^2, \hspace{2ex} z\in\Omega.$ \\

\noindent Taking the supremum over all $f$ in the unit sphere of the Hilbert space ${\mathcal H}_{\Omega}(m_0q\varphi)$ gives the desired estimate. \\

\noindent $(a)$\, We first briefly recall the use made of Skoda's $L^2$ division theorem in [Pop06, Theorem 4.1] to obtain an effective version of a superadditivity result on multiplier ideal sheaves ([Pop06, Theorem 1.2]) corresponding in the present context to the left-hand inclusion in (\ref{eqn:inclusions}). Let $f\in {\mathcal O}(\Omega)$ be an arbitrary element in the unit sphere of the Hilbert space ${\mathcal H}_{\Omega}(m_0(1+\varepsilon)\, \varphi)$. Combined with assumption (\ref{eqn:assumption}), this means that: \\

$\displaystyle 1= \int_{\Omega}|f|^2\, e^{-2m_0(1+\varepsilon)\, \varphi} \, dV_n = \int_{\Omega} \frac{|f|^2}{\bigg(\sum\limits_{j=0}^N|g_j|^2\bigg)^{m_0c(1+\varepsilon)}}\, e^{-2m_0(1+\varepsilon)v}\, dV_n$. \\

\noindent Choose $m_0 \geq \frac{n+2}{c\,\varepsilon}$. We can apply Skoda's $L^2$ division theorem ([Sko72b]) stated as Theorem \ref{The:Skoda} in Appendix \ref{subsection:appendix} to write $f$ as a linear combination with holomorphic coefficients of products of $s:=[m_0c(1+\varepsilon)]-(n+1)$ functions among the $g_j$'s. Namely, for all multi-indices $L=(l_1, \dots , l_s)\in\{1, \dots , N\}^s$, there exist holomorphic functions $h_L$ on $\Omega$ such that:\\

\hspace{6ex}$\displaystyle f=\sum\limits_L h_L g^L \hspace{2ex} \mbox{on} \hspace{2ex} \Omega, \hspace{5ex} \mbox{with} \hspace{2ex} g^L=g_{l_1} \dots g_{l_s},$\\

\noindent with precise $L^2$ estimates on the coefficients $h_L$ (see Theorem \ref{The:Skoda}). Combined with the Cauchy-Schwarz inequality for $|f|^2$ and the submean value inequality for each $|h_L|^2$ on $B\Subset \Omega$, these $L^2$ estimates yield: \\

\hspace{6ex}$|f|^2 \leq C'_{m_0}\,(\sum\limits_{j=0}^N |g_j|^2)^{[m_0c(1+\varepsilon)]-(n+1)}$ \hspace{2ex} on $B$, \\

\noindent with a constant $C'_{m_0}>0$ whose dependence on $m_0$ is explicit. Thus: \\

\hspace{6ex} $|f|^2\, e^{-2m_0\varphi}\leq C'_{m_0}\,(\sum\limits_{j=0}^N |g_j|^2)^{[m_0c(1+\varepsilon)]-(n+1)-m_0c}$ \hspace{2ex} on $B$, \\

\noindent and the crucial fact is that the exponent $[m_0c(1+\varepsilon)]-(n+1)-m_0c$ is non-negative by the choice of $m_0\geq \frac{n+2}{c\,\varepsilon}$. Therefore, the right-hand term above is bounded on $B$ and thus the initial $L^2$ condition satisfied by $f$ on $\Omega$ leads to an $L^{\infty}$ property on $B$ for a slightly less singular weight (i.e. without $(1+\varepsilon)$ in the exponent). The explicit bound we finally get is: \\

\hspace{6ex}$\displaystyle |f|^2\, e^{-2m_0\, \varphi}\leq C_n\,(m_0c(1+\varepsilon)-n)\, (\sup\limits_Be^{\varphi})^{2m_0\varepsilon}:=C_{m_0}$, \\

\noindent on $B\Subset \Omega$, where $C_n>0$ is a constant depending only on $n$ and the diameters of $B$ and $\Omega$. The details can be found in [Pop06] (taking in Theorem 4.1. obtained there $c=1+\varepsilon$, $1-\delta =\frac{1}{1+\varepsilon}$, $m=m_0$). This readily implies that for any $q$ functions $f_1, \dots , f_q$ in the unit sphere of ${\mathcal H}_{\Omega}(m_0(1+\varepsilon)\, \varphi)$ we have: \\

\hspace{6ex} $\displaystyle |f_1 \dots  f_q|^2\, e^{-2m_0q\, \varphi}\leq C_{m_0}^q$ \hspace{3ex} on $B$, \\

\noindent and in particular $f_1 \dots  f_q$ is a section on $B$ of the ideal sheaf ${\mathcal I}(m_0q\, \varphi)$ with $L^2$ norm 

\begin{equation}\label{eqn:Skodaest}\displaystyle \int\limits_B |f_1 \dots  f_q|^2 \, e^{-2m_0q\, \varphi}\, dV_n \leq \mbox{Vol}(B)\, C_{m_0}^q,\end{equation} 

\noindent where $\mbox{Vol}(B)$ stands for the Lebesgue measure of $B$. We will now apply the effective local finite generation theorem for multiplier ideal sheaves (Theorem 1.1 in [Pop06]). That theorem gives the existence of an open subset $B'_0\Subset B$ and of an orthonormal basis $(\sigma_{m_0q, \, j})_{j\in \N}$ of ${\mathcal H}_{\Omega}(m_0q\, \varphi)$ such that finitely many elements of this basis, say $\sigma_{m_0q, \, j},$ $1\leq j\leq N_{m_0q},$ generate ${\mathcal I}(m_0q\, \varphi)$ on $B'_0$ in such a way that estimates analogous to (\ref{eqn:redecomp})--(\ref{eqn:recomparison}) hold. When applied to $f_1 \dots  f_q$ regarded as a section on $B$ of ${\mathcal I}(m_0q\varphi)$ with local generators $\sigma_{m_0q, \, j}, 1\leq j\leq N_{m_0q}$, (\ref{eqn:redecomp}) reads:

\begin{equation}\label{eqn:generation}\displaystyle f_1(z) \dots  f_q(z)= \sum\limits_{j=1}^{N_{m_0q}}h_{j_1, \dots , j_q}^{(j)}(z) \, \sigma_{m_0q, \, j}(z), \hspace{3ex} z\in B'_0, \,\, m_0, q \gg  1\end{equation}

\noindent with holomorphic coefficients $h_{j_1, \dots , j_q}^{(j)}$ estimated as (cf.(\ref{eqn:reestimate}), (\ref{eqn:Skodaest})): 

\begin{eqnarray}\nonumber\displaystyle \sup\limits_{B'_0}\sum\limits_{j=1}^{N_{m_0q}}|h_{j_1, \dots , j_q}^{(j)}|^2 & \leq & \displaystyle C\, N_{m_0q}\, \int\limits_B |f_1 \dots  f_q|^2 \, e^{-2m_0q\, \varphi} \\
 \label{eqn:hcoeff}& \leq & C\, N_{m_0q}\,\mbox{Vol}(B)\, C_{m_0}^q, \hspace{3ex} m_0, q \gg 1.
\end{eqnarray}

\noindent Taking derivatives in identity (\ref{eqn:generation}) we get:

\begin{equation}\label{eqn:derivgeneration}\displaystyle D^{\alpha}(f_1 \dots  f_q)= \sum\limits_{j=1}^{N_{m_0q}}\sum\limits_{\beta\leq \alpha}{\alpha \choose \beta} D^{\alpha-\beta}h_{j_1, \dots , j_q}^{(j)}\, D^{\beta} \sigma_{m_0q, \, j} \hspace{2ex} \mbox{on} \,\, B'_0.\end{equation}

\noindent Applying twice the Cauchy-Schwarz inequality we get at every point in $B'_0$:

\begin{equation}\label{eqn:derivineq}\displaystyle |D^{\alpha}(f_1 \dots  f_q)|^2 \leq C'_{\alpha}\, \sum\limits_{\beta\leq \alpha}\bigg(\sum\limits_{j=1}^{N_{m_0q}}|D^{\alpha-\beta}h_{j_1, \dots , j_q}^{(j)}|^2\bigg)\, \bigg(\sum\limits_{j=1}^{N_{m_0q}}|D^{\beta}\sigma_{m_0q, \, j}|^2\bigg),\end{equation}

\noindent where $C'_{\alpha}:=\sum\limits_{\beta\leq \alpha}{\alpha \choose \beta}^2$. The supremum of the term on the left in the above estimate (\ref{eqn:derivineq}) over all $f_1, \dots , f_q$ ranging over the unit sphere of ${\mathcal H}_{\Omega}(m_0(1+\varepsilon)\, \varphi)$ equals at every point in $\Omega$ the value of the left-hand term in the inequality we intend to prove. The estimate of the coefficients $h_{j_1, \dots , j_q}^{(j)}$ obtained in (\ref{eqn:hcoeff}), combined with the Cauchy inequalities estimating $|D^{\beta}h_{j_1, \dots , j_q}^{(j)}|^2$ above on any $B_0\Subset B'_0$ in terms of $\sup\limits_{B'_0}|h_{j_1, \dots , j_q}^{(j)}|^2$, gives the desired estimate on $B_0$ . \hfill $\Box$

\vspace{3ex}

The main ingredient in the regularisation process to be described in the next section will be the following estimate on derivatives relying on Theorem \ref{The:potentiel} combined with the Ohsawa-Takegoshi $L^2$ extension theorem applied on a complex line. This strategy was already used back in section \ref{subsection:local} when the Lelong numbers of $\varphi$ were assumed to vanish. The extra difficulty in the general case stems from the fact that in Theorem \ref{The:potentiel} the distances between the points $a_j$ where the Lelong numbers of $\varphi$ are $\geq \frac{1}{m}$ cannot be estimated when $m\rightarrow +\infty$ (they may decrease arbitrarily fast to $0$). As the derivatives of the functions $f_m$ constructed in Theorem \ref{The:potentiel} depend on these distances, it follows that we cannot control the growth of these derivatives as $m\rightarrow +\infty$ if this procedure is applied to $m\varphi$ to produce sections of ${\mathcal I}(m\varphi)$. The solution we propose to this problem is to choose $m=m_0q$ and to apply this procedure to  $m_0(1+\varepsilon)\varphi$ instead in order to construct sections $g_{m_0(1+\varepsilon)}$ of ${\mathcal I}(m_0(1+\varepsilon)\varphi)$. Then the left-hand inclusion in (\ref{eqn:inclusions}) of Proposition \ref{Prop:inclusions} ensures that $g_{m_0(1+\varepsilon)}^q$ is a section of ${\mathcal I}(m_0q\varphi)={\mathcal I}(m\varphi)$ and we have a complete control on it and its derivatives by means of the effective estimates of Proposition \ref{Prop:effectiveincl}. The points $a_j$ featuring in the definition of $g_{m_0(1+\varepsilon)}^q$ are now the same as those of $g_{m_0(1+\varepsilon)}$ (only the exponents $m_j$ are multiplied by $q$, see notation in Theorem \ref{The:potentiel}). Thus their minimum mutual distance $\delta_{m_0}$ depends only on $m_0$ (though in an uncontrollable fashion). To obtain a control of the Monge-Amp\`ere masses in the next section, we shall choose $m=m_0q$ with $q=q(m_0)$ sufficiently large to neutralise the growth of $\delta_{m_0}\downarrow 0$. The main interest of the next proposition is that it gives a lower bound independent of $q$.

\begin{Prop}\label{Prop:derivlowerbound} Let $\varphi$ be any psh function on $\Omega\Subset \C^n$ and let $\nu(x)=\nu(\varphi,\, x)$ denote the Lelong number of $\varphi$ at any $x\in\Omega$. Then, for every $m_0, q\in \N$ and every $0<\varepsilon<1$, any orthonormal basis $(\sigma_{m_0(1+\varepsilon), \, j})_{j\in \N}$ of ${\mathcal H}_{\Omega}(m_0(1+\varepsilon)\varphi)$ satisfies the following estimate at every $x\in \Omega$: \\

 $\displaystyle \frac{1}{2m_0q}\, \log\sum\limits_{j_1, \dots , j_q=0}^{+\infty}\sum\limits_{|\alpha|=0}^{d_{m_0q}(x)}|D^{\alpha}(\sigma_{m_0(1+\varepsilon), \, j_1} \dots \sigma_{m_0(1+\varepsilon), \, j_q})(x)|^2 \geq C_0\, \log \delta_{m_0},$ \\

\noindent where $d_{m_0q}(x):=\max\{[m_0q\nu(x)(1+\varepsilon)], \, 1\}$ and $0<C_0\leq C_{\Omega}\,\int_{\Omega}dd^c\varphi\wedge \beta^{n-1}$ with $C_{\Omega}>0$ depending only on $\Omega$. The constant $\delta_{m_0}>0$ depends only on $m_0$ in a way that we cannot control.

\end{Prop}

\noindent {\it Proof.} Let $x\in \Omega$ be an arbitrary point. The case $\nu(x)=0$ was settled in section \ref{subsection:local}. Assume that $\nu(x)>0$. For every multi-index $\alpha=(\alpha_1, \dots , \alpha_n)\in \N^n$, the linear map of evaluation at $x$:  \\

$\displaystyle {\mathcal H}_{\Omega}(m_0(1+\varepsilon)\varphi)\hat{\otimes}\dots \hat{\otimes} {\mathcal H}_{\Omega}(m_0(1+\varepsilon)\varphi)\ni u \longmapsto D^{\alpha}u(x)\in \C$ \\

\noindent of the derivatives of elements $u$ in the completed tensor product of the Hilbert space ${\mathcal H}_{\Omega}(m_0(1+\varepsilon)\varphi)$ by itself $q$ times defines a continuous linear map whose squared norm can be expressed in two ways as: \\

$\displaystyle\sum\limits_{j_1, \dots , j_q=0}^{+\infty} |D^{\alpha}(\sigma_{m_0(1+\varepsilon), \, j_1} \dots \sigma_{m_0(1+\varepsilon), \, j_q})(x)|^2=\sup\limits_u|D^{\alpha}u(x)|^2$,\\

\noindent since $(\sigma_{m_0(1+\varepsilon), \, j_1} \dots \sigma_{m_0(1+\varepsilon), \, j_q})_{j_1, \dots , j_q\in \N}$ defines an orthonormal basis in the completed tensor product Hilbert space, and the supremum on the right is taken over all $u$ in the unit ball of this space. To find the desired lower bound for the term on the right side, we shall produce an element in the unit ball of the completed tensor product Hilbert space ${\mathcal H}_{\Omega}(m_0(1+\varepsilon)\varphi)\hat{\otimes}\dots \hat{\otimes} {\mathcal H}_{\Omega}(m_0(1+\varepsilon)\varphi)$ ($q$ times) for which one of the partial derivatives up to order $[A\, m_0(1+\varepsilon)q]$ can be estimated below in absolute value at $x$. This element will be chosen as the $q^{th}$ power of an element in the unit ball of ${\mathcal H}_{\Omega}(m_0(1+\varepsilon)\varphi)$.

 Let $L$ be a complex line through $x$ such that the restriction of $\varphi$ to $L$ has the same Lelong number at $x$ as $\varphi$. This is the case for almost all lines passing through $x$ ([Siu74]). After possibly changing coordinates, we can assume that $x=0$ and $L=\{z_2= \dots = z_n=0\}$. As in the Introduction, we have a decomposition: \\

\hspace{6ex} $\displaystyle \varphi_{|L} = N\star \Delta \varphi_{|L} + \mbox{Re}\, g$  \hspace{3ex} on $\Omega \cap L$, \\

\noindent for some holomorphic function $g$ and the one-dimensional Newton kernel $N$. Theorem \ref{The:potentiel} applied to $(1+\varepsilon)\, \varphi_{|L}$ gives the existence, for every $m_0\in \N$, of a holomorphic function of one variable $g_{m_0(1+\varepsilon)}$ on $\Omega \cap L$ of the form: \\

\hspace{6ex} $\displaystyle g_{m_0(1+\varepsilon)}(z_1) = e^{m_0(1+\varepsilon)\, g(z_1)}\, \prod\limits_{j=1}^{N_{m_0}}(z-a_j)^{m_j}$, \hspace{3ex} $\displaystyle z_1\in \Omega \cap L$, \\

\noindent such that $\sum\limits_{j=1}^{N_{m_0}}m_j\leq C_0\, (1+\varepsilon)m_0$ with a constant $C_0=\int\limits_{\Omega\cap L}dd^c\varphi_{|L}>0$ independent of $m_0$. It further satisfies: \\

\hspace{6ex} $\displaystyle C_{m_0}:=\int_{\Omega\cap L}|g_{m_0(1+\varepsilon)}|^2\, e^{-2m_0(1+\varepsilon)\varphi}\, dV_L = o(m_0),$ \\

\noindent where $dV_L$ is the volume form on $L$. We can now apply the Ohsawa-Takegoshi $L^2$ extension theorem ([Ohs88, Corollary 2, p. 266]) to get a holomorphic extension $G_{m_0(1+\varepsilon)}\in {\mathcal H}_{\Omega}(m_0(1+\varepsilon)\varphi)$ of $g_{m_0(1+\varepsilon)}$ from $\Omega\cap L$ to $\Omega$ which satisfies the $L^2$ estimate: \\

$\displaystyle \int\limits_{\Omega}|G_{m_0(1+\varepsilon)}|^2\, e^{-2m_0(1+\varepsilon)\varphi}\, dV_n \leq C_n\, \int\limits_{\Omega\cap L}|g_{m_0(1+\varepsilon)}|^2\, e^{-2m_0(1+\varepsilon)\varphi}\, dV_L=C_n\, C_{m_0},$ \\

\noindent with a constant $C_n>0$ depending only on $n$. Thus, $\frac{1}{(C\, C_{m_0})^{1/2}}\, G_{m_0(1+\varepsilon)}$ belongs to the unit ball of ${\mathcal H}_{\Omega}(m_0(1+\varepsilon)\varphi)$. Then the holomorphic function defined on $\Omega$ as: \\

\hspace{10ex} $\displaystyle F_{m_0q}:= \frac{1}{(C_n\, C_{m_0})^{q/2}} \, G_{m_0(1+\varepsilon)}^q$ \\

\noindent belongs to the unit ball of the completed tensor product Hilbert space ${\mathcal H}_{\Omega}(m_0(1+\varepsilon)\varphi)\hat{\otimes}\dots \hat{\otimes} {\mathcal H}_{\Omega}(m_0(1+\varepsilon)\varphi)$ ($q$ times). Suppose, for example, that $x=a_k$ for some $k$ and let $\delta_{m_0}:=\min\limits_{j\neq k}|a_j-a_k|>0$. As, by construction, the derivatives  $D^{(l, \, 0 \dots ,\, 0)} G_{m_0(1+\varepsilon)}(x)$ in the direction of the line $L$ coincide on $L$ with the derivatives $g_{m_0(1+\varepsilon)}^{(l)}(x)$, $l\in\N$, we get, for an appropriate $l\leq d_{m_0q}(x)$: \\

\hspace{6ex}$\bigg|F_{m_0q}^{(l)}(x)\bigg|\geq \frac{e^{qm_0(1+\varepsilon)\, Re\, g(z_1)}}{(C_n\, C_{m_0})^{q/2}} \,\, \delta_{m_0}^{q\,\big(\sum_{j=1}^{N_{m_0}}m_j\big)} \geq \frac{e^{qm_0(1+\varepsilon)\, Re\, g(z_1)}}{(C_n\, C_{m_0})^{q/2}} \,\, \delta_{m_0}^{C_0(1+\varepsilon)qm_0}.$ \\

\noindent Then $\frac{1}{2m_0q}\, \log |F_{m_0q}^{(l)}(x)|\geq C_0(1+\varepsilon)\, \log\delta_{m_0}+ O(1)$, and the stated lower bound follows after we absorb $(1+\varepsilon)$ in $C_0$ and the $O(1)$ term in the first term. The estimate of the constant $C_0$ follows from the following averaging argument. Suppose, for simplicity, that $\Omega=B(0,\,r)$ for some $r>0$. Then:\\

$\begin{array}{lll}\int\limits_{L\in \Proj^{n-1}}(\int\limits_{B(0,\,r)\cap L}dd^c\varphi_{|L})\, dv(L) & = & \int\limits_{B(0,\,r)}dd^c\varphi\wedge\int\limits_{L\in\Proj^{n-1}}[L]\, dv(L) \\

  =\int\limits_{B(0,\,r)}dd^c\varphi\wedge (dd^c\log |z|)^{n-1} & = & \frac{1}{r^{2(n-1)}}\, \int\limits_{B(0,\,r)}dd^c\varphi\wedge \beta^{n-1},\end{array},$

\vspace{2ex}

\noindent where $[L]$ denotes the $(n-1,\,n-1)$-current of integration on the line $L$, $dv$ denotes the unique unitary invariant measure of total mass $1$ on projective space, the second equality above follows from Crofton's formula (see e.g. [Dem97, chapter III, p.196]), and the third equality is a well-known formula of Lelong (see e.g. [Dem97, chapter III, Formula 5.5.]). In particular, for $L$ in a subset of $\Proj^{n-1}$ of positive $dv$-measure, $C_0= \int\limits_{B(0,\,r)\cap L}dd^c\varphi_{|L} \leq \frac{2}{r^{2(n-1)}}\, \int\limits_{B(0,\,r)}dd^c\varphi\wedge \beta^{n-1}$, and it is enough to choose such a line $L$ to have the stated estimate on $C_0$. \hfill $\Box$ \\

The estimates obtained in Proposition \ref{Prop:effectiveincl} $(a)$ and in Proposition \ref{Prop:derivlowerbound} combine to achieve the main purpose of this section: a lower estimate for derivatives in the following sense. Given pseudoconvex open sets $B\Subset \Omega\Subset \C^n$ and a psh function $\varphi$ on $\Omega$ of the form (\ref{eqn:assumption}), let $\nu:=\sup\limits_{x\in B}\nu(\varphi,\, x)$. Suppose that $\nu>0$. If $(\sigma_{m,\, j})_{j\in\N}$ is an orthonormal basis of ${\mathcal H}_{\Omega}(m\varphi)$ such that finitely many elements $(\sigma_{m,\, j})_{1\leq j\leq N_{m}}$ generate ${\mathcal I}(m\varphi)$ effectively on some $B_0\Subset B$ with properties (\ref{eqn:redecomp})---(\ref{eqn:recomparison}) satisfied, we obtain a lower estimate on $B_0$ for the psh functions:

\begin{equation}\label{eqn:regderiv}\displaystyle u_m:=\frac{1}{2m}\, \log\sum\limits_{j=1}^{N_m}\sum\limits_{|\alpha|=0}^{[m\nu(1+\varepsilon)]} \bigg|\frac{D^{\alpha}\sigma_{m,\,j}}{\alpha!}\bigg|^2\end{equation}

\noindent as $m\rightarrow +\infty$ in the following form summing up the discussion in this section.

\begin{Cor}\label{Cor:below} For all $q\in \N$, $0<\varepsilon << 1,$ and $m_0\geq \frac{n+2}{c\,\varepsilon}$, there exists an orthonormal basis $(\sigma_{m_0q,\, j})_{j\in\N}$ of ${\mathcal H}_{\Omega}(m_0q\varphi)$ satisfying, for $m=m_0q$, the estimate: \\

$\displaystyle u_m:=\frac{1}{2m}\, \log\sum\limits_{j=1}^{N_{m}}\sum\limits_{|\alpha|=0}^{[m\nu(1+\varepsilon)]}\bigg|\frac{D^{\alpha}\sigma_{m,\,j}}{\alpha!}\bigg|^2 \geq  C_0 \log\delta_{m_0} - A_m$ \hspace{2ex} on \,\, $B_0.$ \\

\noindent Crucially, $\delta_{m_0}$ depends only on $m_0$ (though in an uncontrollable way) and not directly on $m$. The other constants satisfy: $0\leq C_0\leq C_{\Omega}\, \int\limits_{\Omega}dd^c\varphi\wedge \beta^{n-1},$ $0\leq A_m\leq \frac{1}{m}\, \log(C_n\,(m\nu)^n\, N_m^2 \, C_{m_0}^q)$, with $C_{m_0}=C_n\, (m_0c(1+\varepsilon)-n)(\sup\limits_B e^{\varphi})^{2m_0\varepsilon},$ and $C_n$ depending only on $n,\, B_0,\, B,\, \Omega$.

\end{Cor}

\noindent {\it Proof.} This is done in two steps. First, apply Proposition \ref{Prop:effectiveincl} $(a)$ for every $\alpha\in\N^n$, sum over $0\leq |\alpha|\leq [m(1+\varepsilon)\nu]$, and take $\frac{1}{2m}\log$ to get the following lower estimate on $u_m=u_{m_0q}$ for any orthonormal basis $(\sigma_{m_0(1+\varepsilon),\, j})_{j\in\N}$ of ${\mathcal H}_{\Omega}(m_0(1+\varepsilon)\varphi):$ \\

 $\displaystyle u_m \geq \frac{1}{2m}\, \log\sum\limits_{j_1, \dots , j_q=0}^{+\infty}\sum\limits_{|\alpha|=0}^{[m_0q\nu(1+\varepsilon)]}\bigg|\frac{D^{\alpha}(\sigma_{m_0(1+\varepsilon),\, j_1}\dots \sigma_{m_0(1+\varepsilon),\, j_q})}{\alpha!}\bigg|^2 - A_m$ \\

\noindent on $B_0$. To estimate $A_m$ we notice, with the notation of Proposition \ref{Prop:effectiveincl}$(a)$, that: $\sum\limits_{|\alpha|=0}^{[m(1+\varepsilon)\nu]}C_{\alpha}\leq C_n'\, m^n(1+\varepsilon)^n\nu^n\leq 2^nC_n'\, (m\nu)^n$, for some $C'_n>0$. Second, apply Proposition \ref{Prop:derivlowerbound} to get the desired lower bound for the term appearing on the right above. \hfill $\Box$ \\

 The above functions $u_m$ modify Demailly's regularisation of $\varphi$ (cf. (\ref{eqn:Demaillyapprox})) by taking into account derivatives up to order $[m(1+\varepsilon)\nu]$ and only finitely many $\sigma_{m,\,j}$'s. The upshot is that the functions $u_{m}$ can be estimated below by finite constants. This is a major improvement of the situation described in Lemma \ref{Lem:newreg} where the lower bound for $\psi_m$ in terms of $\varphi$ was unsatisfactory owing to the $-\infty$ poles of $\varphi$. The present functions $u_m$ are thus candidates to defining better regularisations of $\varphi$. However, the shortcoming of such a definition would be that $u_m$ may not converge to $\varphi$ as $m\rightarrow +\infty$ owing to the derivation order $[m(1+\varepsilon)\nu]$ being too large (the limit of $[m(1+\varepsilon)\nu]/m$ is $>0$ instead of zero as required by Lemma \ref{Lem:derivorder} for convergence).

We will overcome this obstacle in the next section by applying this procedure after we have removed enough singularities to lower substantially the necessary derivation order.

\section{Regularisation with mass control: the general case}\label{subsection:reg}

Let $\varphi$ be an arbitrary psh function on a bounded pseudoconvex domain $\Omega\Subset \C^n$ such that $i\partial\bar{\partial}\varphi \geq C_0\, \beta$ for some $C_0>0$. As usual, $\beta$ is the standard K\"ahler form on $\C^n$. For every fixed $p\in \N$, we will approximate Demailly's regularising function $\varphi_p$ (cf. (\ref{eqn:Demaillyapprox})): \\

 $\displaystyle \varphi_p=\frac{1}{2p}\, \log\sum\limits_{j=0}^{+\infty}|\sigma_{p,\,j}|^2$, \hspace{2ex} $(\sigma_{p,\,j})_{j\in\N}$ an orthonormal basis of ${\mathcal H}_{\Omega}(p\varphi),$\\ 

 \noindent by a sequence of psh functions $(\psi_{m, \, p})_{m\in \N}$ with analytic singularities for which we can control the Monge-Amp\`ere masses. In the last step of the proof, the desired approximation of the original $\varphi$ will be obtained by letting $p\rightarrow +\infty$. The results of the previous section \ref{subsection:addefect} will now be applied to $\varphi_p$ (which is of the form (\ref{eqn:assumption}) with $c=\frac{1}{p}$) in place of $\varphi$. The notation is the same as in section \ref{subsection:addefect} with an extra $p$ in the indices (e.g. $(\sigma_{m_0(1+\varepsilon),\,p,\,j})_{j\in\N}$ denotes an orthonormal basis of ${\mathcal H}_{\Omega}(m_0(1+\varepsilon)\varphi_p)$, $N_{m_0q,\, p}$ is the number of local generators of ${\mathcal I}(m_0q\varphi_p)$, etc).

 The $m^{th}$ regularisation $\psi_{m,\,p}$ of $\varphi_p$ will be constructed using the ideal sheaf ${\mathcal I}(m\varphi_p)$. As already explained, the main idea is to write $m=m_0q$ with $q\gg m_0$ and to make full use of the inclusions (\ref{eqn:inclusions}) of Proposition \ref{Prop:inclusions} which show intuitively that ${\mathcal I}(m\varphi_p)$ and ${\mathcal I}(m_0\varphi_p)^q$ are ``almost'' equal when $m_0 \rightarrow +\infty$ (and $\varepsilon \downarrow 0$ with $m_0\geq p\,\frac{n+2}{\varepsilon}$). Considering a $\log$-resolution of the ideal sheaf ${\mathcal I}(m_0\varphi_p)$, we shall ``almost'' get a $\log$-resolution of ${\mathcal I}(m\varphi_p)$ with $m=m_0\, q \gg  m_0$ up to a small loss which can be controlled explicitly in terms of $0<\varepsilon << 1$ by means of the estimates of Propositions \ref{Prop:effectiveincl} and \ref{Prop:derivlowerbound}. The advantage of this method over considering right away a $\log$-resolution of ${\mathcal I}(m\varphi_p)$ is that the blow-up does not depend explicitly on $m$ but only on $m_0$. Thus the mass estimates of the $m^{th}$ regularisation $\psi_{m,\,p}$ will only depend on $m_0$ (and $p$) and can be neutralised by choosing $m=m_0q$ with $q=q(m_0)$ large enough. This idea ties in with the explanation given before Proposition \ref{Prop:derivlowerbound}. \\

 Let $\tilde{\Omega}=\tilde{\Omega}_{m_0}$ be a smooth variety and let $\mu=\mu_{m_0}:\tilde{\Omega}\rightarrow \Omega$ be a proper modification (i.e. a holomorphic bimeromorphic map) such that \\

\hspace{3ex}$\displaystyle \mu^{\star}{\mathcal I}(m_0(1+\varepsilon)\varphi_p)={\mathcal O}(-E_{m_0(1+\varepsilon),\, p})$, \hspace{3ex} $\displaystyle \mu^{\star}{\mathcal I}(m_0\varphi_p)={\mathcal O}(-E_{m_0,\, p}),$ \\

\noindent where $E_{m_0(1+\varepsilon),\, p}, E_{m_0,\, p}$ are effective normal crossing divisors on $\tilde{\Omega}$. If $V\,{\mathcal I}(m_0(1+\varepsilon)\varphi_p)$ denotes the zero variety of ${\mathcal I}(m_0(1+\varepsilon)\varphi_p)$, the restriction \\

\hspace{6ex} $\displaystyle\mu:\tilde{\Omega}\setminus E_{m_0(1+\varepsilon),\, p} \longrightarrow \Omega\setminus V\,{\mathcal I}(m_0(1+\varepsilon)\varphi_p)$ \\

\noindent defines a biholomorphism. The inclusion $\mu^{\star}{\mathcal I}(m_0(1+\varepsilon)\varphi_p)\subset \mu^{\star}{\mathcal I}(m_0\varphi_p)$ also implies the existence of an effective divisor $D_{m_0(1+\varepsilon),\, p}$ such that \\

\hspace{6ex} $E_{m_0(1+\varepsilon),\, p}=E_{m_0,\, p} + D_{m_0(1+\varepsilon),\, p}$. \\

\noindent Let $(\tilde{U}_l)_{1\leq l \leq N}$ be a finite collection of open coordinate balls covering $\tilde{\Omega}$ such that the restrictions of $\mu^{\star}{\mathcal I}(m_0(1+\varepsilon)\varphi_p)$ and $\mu^{\star}{\mathcal I}(m_0\varphi_p)$ to each of these have a unique generator. Let $(U_l)_{1\leq l \leq N}$ be a corresponding collection of open balls covering $\Omega$ such that $\tilde{U}_l \setminus E_{m_0(1+\varepsilon),\, p}=\mu^{-1}(U_l\setminus V\,{\mathcal I}(m_0(1+\varepsilon)\varphi_p))$. The modification $\mu$ is a simultaneous $\log$-resolution of the ideal sheaves ${\mathcal I}(m_0(1+\varepsilon)\varphi_p)$ and ${\mathcal I}(m_0\varphi_p)$ whose existence follows from Hironaka's theorem on the resolution of singularities.

 In the sequel, we shall need the following comparison lemma in passing from the global to the local picture on coordinate patches and vice-versa.

\begin{Lem}\label{Lem:locglobal} Let $\varphi$ be a psh function on a bounded pseudoconvex open set $\Omega \subset \C^n$, and let $B\Subset \Omega$ be a relatively compact open subset. For every $m,p\in \N$, the expressions analogous to the Bergman kernels associated with the weight $m\varphi$ on $\Omega$ and respectively $B$ which take into account derivatives up to order $p$: \\

\hspace{6ex} $\displaystyle B^{(p)}_{m\varphi, \, \Omega}:=\sum\limits_{k=0}^{+\infty}\sum\limits_{|\alpha|=0}^p|D^{\alpha}\sigma_{m, \, k}|^2,$ \hspace{3ex} $\displaystyle B^{(p)}_{m\varphi, \, B}:=\sum\limits_{k=0}^{+\infty}\sum\limits_{|\alpha|=0}^p|D^{\alpha}\mu_{m, \, k}|^2,$ \\

\noindent defined by orthonormal bases $(\sigma_{m, \, k})_{k\in \N}$ and $(\mu_{m, \, k})_{k\in \N}$ of the Hilbert spaces ${\mathcal H}_{\Omega}(m\varphi)$ and respectively ${\mathcal H}_{B}(m\varphi)$, can be compared as: \\

\hspace{3ex} $\displaystyle B^{(p)}_{m\varphi, \, \Omega} \leq B^{(p)}_{m\varphi, \, B} \leq C_{n, \, d, \, r}\, (d/r)^p\, B^{(p)}_{m\varphi, \, \Omega}$ \hspace{2ex} on any $B_0\Subset B\Subset \Omega,$ \\

\noindent where $C_{n, \, d, \, r}>0$ is a constant depending only on $n$, the diameter $d$ of $\Omega$, and the distance $r>0$ between the boundaries of $B_0$ and $B$. \\

\end{Lem}

\noindent {\it Proof.} It is by a standard application of H\"ormander's $L^2$ estimates ([Hor65]) and it is given in Appendix \ref{subsection:appendix}. \hfill $\Box$

\vspace{3ex}

\noindent We will now concentrate attention on an arbitrary $\tilde{U}_l$ that we generically call $\tilde{U}$. Let $g_{m_0(1+\varepsilon),\, p}$ (respectively $g_{m_0,\, p}$) be the unique generator on $\tilde{U}$ of $\mu^{\star}{\mathcal I}(m_0(1+\varepsilon)\varphi_p)$ (respectively $\mu^{\star}{\mathcal I}(m_0\varphi_p)$). By pull-back, the inclusions (\ref{eqn:inclusions}) become:

\begin{equation}\mu^{\star}{\mathcal I}(m_0(1+\varepsilon)\varphi_p)_{|\tilde{U}}^q \subset \mu^{\star}{\mathcal I}(m_0q \varphi_p)_{|\tilde{U}} \subset \mu^{\star}{\mathcal I}(m_0 \varphi_p)_{|\tilde{U}}^q.\end{equation}

\noindent The analogous inclusions to (\ref{eqn:inclusions}) of Proposition \ref{Prop:inclusions} applied to $\varphi_p\circ \mu$ on $\tilde{U}\Subset \tilde{U}'$ in place of $\varphi_p$ on $B\Subset \Omega$ read:

\begin{equation}\label{eqn:upinclusions}{\mathcal I}(m_0(1+\varepsilon)\varphi_p\circ \mu)_{|\tilde{U}}^q \subset {\mathcal I}(m_0q \varphi_p\circ \mu)_{|\tilde{U}} \subset {\mathcal I}(m_0 \varphi_p\circ \mu)_{|\tilde{U}}^q,\end{equation}

\noindent with the three ideal sheaves in (\ref{eqn:upinclusions}) above generated respectively on $\tilde{U}$ by: \\

\noindent$\begin{array}{lll}\displaystyle(\tilde{\sigma}_{m_0(1+\varepsilon), \, p, \, j_1} \dots \tilde{\sigma}_{m_0(1+\varepsilon), \, p, \, j_q})_{j_1, \dots , j_q\in \N}, & \mbox{with} & \displaystyle \tilde{\sigma}_{m_0(1+\varepsilon), \, p, \, j}:=J_{\mu} \, \sigma_{m_0(1+\varepsilon),\, p, \, j}\circ \mu \\

\displaystyle (\tilde{\sigma}_{m_0q, \, p, \, j})_{j\in \N}, & \mbox{with} & \displaystyle \tilde{\sigma}_{m_0q, \, p, \, j}:=J_{\mu} \, \sigma_{m_0q,\,p,\,j}\circ \mu, \\

\displaystyle (\tilde{\sigma}_{m_0, \, p, \, j_1} \dots \tilde{\sigma}_{m_0, \, p, \, j_q})_{j_1, \dots , j_q\in \N}, & \mbox{with} & \displaystyle \tilde{\sigma}_{m_0, \, p, \, j}:=J_{\mu} \, \sigma_{m_0, \, p, \, j}\circ \mu,\end{array}$

\vspace{2ex}

\noindent where $J_{\mu}$ is the Jacobian of $\mu$. On the other hand, the ideal sheaves ${\mathcal I}(m_0(1+\varepsilon)\varphi_p\circ \mu)$ and ${\mathcal I}(m_0\varphi_p\circ \mu)$ have on $\tilde{U}$ unique generators: \\

\hspace{6ex}$\displaystyle\tilde{g}_{m_0(1+\varepsilon), \, p}:=J_{\mu}\, g_{m_0(1+\varepsilon), \, p}$ \hspace{2ex} and respectively \hspace{2ex}  $\displaystyle\tilde{g}_{m_0, \, p}:=J_{\mu}\, g_{m_0, \, p}$. \\

\noindent Thus, on $\tilde{U}$, we get: \\

\hspace{6ex} $\displaystyle \tilde{\sigma}_{m_0(1+\varepsilon),\, p,\, j}=\tilde{g}_{m_0(1+\varepsilon),\, p}\, \tilde{\mu}_{m_0(1+\varepsilon),\, p,\, j},$ \hspace{3ex}  $\displaystyle \tilde{\sigma}_{m_0,\, p,\, j}=\tilde{g}_{m_0,\, p}\, \tilde{\mu}_{m_0,\, p,\, j},$ \\

\noindent with holomorphic functions $(\tilde{\mu}_{m_0(1+\varepsilon),\, p,\, j})_{j\in \N}$ without common zeroes and holomorphic functions $(\tilde{\mu}_{m_0,\, p,\, j})_{j\in \N}$ without common zeroes as well. Moreover, by the right-hand inclusion in (\ref{eqn:upinclusions}) there are holomorphic functions $(\tilde{\mu}_{m_0q,\, p,\, j})_{j\in \N}$ on $\tilde{U}$ such that: \\

\hspace{6ex} $\displaystyle \tilde{\sigma}_{m_0q,\, p,\, j}=\tilde{g}_{m_0,\, p}^q \, \tilde{\mu}_{m_0q,\, p,\, j}$, \hspace{3ex} $j\in \N$. \\

\noindent The functions $(\tilde{\mu}_{m_0q,\, p,\, j})_{j\in \N}$ may still have common zeroes to the extent to which the right-hand inclusion in (\ref{eqn:upinclusions}) fails to be an equality. However, the conjunction with the left-hand inclusion in (\ref{eqn:upinclusions}) limits the amount of common zeroes which we are now proceeding to estimate.

 By virtue of Lemma \ref{Lem:locglobal}, replacing ${\mathcal H}_{\Omega}(m_0q\varphi_p)$ with its counterpart ${\mathcal H}_U(m_0q\varphi_p)$ defined on a smaller set $U\Subset \Omega$ can only change the Bergman kernel estimates and their analogues with derivatives by insignificant constants. We can thus suppose, without loss of generality, that $(\sigma_{m_0q, \, p, \, j})_{j\in \N}$ is an orthonormal basis of ${\mathcal H}_U(m_0q\varphi_p)$ (and the analogous assumptions for the weights $m_0(1+\varepsilon)\varphi_p$ and $m_0\varphi_p$.) The change of variable formula shows that: \\

$\begin{array}{lll} 1 & = & \displaystyle\int\limits_{U}|\sigma_{m_0q, \, p, \, j}|^2\, e^{-2m_0q\varphi_p}\, dV=\int\limits_{\tilde{U}}|J_{\mu}|^2 \, |\sigma_{m_0q, \, p, \, j}\circ \mu|^2\, e^{-2m_0q\varphi_p\circ \mu}\, d\tilde{V} \\

  & = & \displaystyle\int\limits_{\tilde{U}} |\tilde{\mu}_{m_0q, \, p, \, j}|^2\, e^{-2m_0q (\varphi_p\circ \mu - \frac{1}{m_0}\, \log|\tilde{g}_{m_0, \, p}|)}\, d\tilde{V}.\end{array}$

\vspace{1ex}

\noindent If we denote: \hspace{3ex} $\displaystyle \tilde{\psi}_{m_0, \, p}:=\varphi_p\circ \mu - \frac{1}{m_0}\, \log|\tilde{g}_{m_0, \, p}|$ \hspace{2ex} on $\tilde{U}$, \\

\noindent we see that $\tilde{\psi}_{m_0, \, p}$ is a psh function on $\tilde{U}$ and $(\tilde{\mu}_{m_0q, \, p, \, j})_{j\in \N}$ defines an orthonormal basis of ${\mathcal H}_{\tilde{U}}(m_0q\, \tilde{\psi}_{m_0, \, p}).$ The same argument shows that $(\tilde{\sigma}_{m_0, \, p, \, j})_{j\in \N}$ defines an orthonormal basis of ${\mathcal H}_{\tilde{U}}(m_0\varphi_p\circ \mu)$ and implicitly a system of generators for the ideal sheaf ${\mathcal I}(m_0\varphi_p\circ \mu)$ on $\tilde{U}$. Demailly's inequality on Lelong numbers (\ref{eqn:Lelongnumbers}) applied to ${\mathcal H}_{\tilde{U}}(m_0\varphi_p\circ \mu)$ shows that the unique generator $\tilde{g}_{m_0, \, p}$ of ${\mathcal I}(m_0\varphi_p\circ \mu)$ (which concentrates all the common zeroes of the $\tilde{\sigma}_{m_0, \, p, \, j}$'s) satisfies: \\

$\displaystyle \nu(\varphi_p\circ \mu, \, x) -\frac{n}{m_0} \leq \frac{1}{m_0}\, \nu(\log|\tilde{g}_{m_0, \, p}|, \, x)\leq \nu(\varphi_p\circ \mu, \, x), \hspace{2ex} x\in \tilde{U},$ \\

\noindent and implicitly \\

\hspace{6ex} $\displaystyle 0\leq \nu(\tilde{\psi}_{m_0, \, p}, \, x) \leq \frac{n}{m_0}, \hspace{2ex} x\in \tilde{U}.$ \\

\noindent By the same inequality (\ref{eqn:Lelongnumbers}) applied to ${\mathcal H}_{\tilde{U}}(m_0q\, \tilde{\psi}_{m_0, \, p})$, the minimum of the vanishing orders of the $\tilde{\mu}_{m_0q, \, p, \, j}$'s satisfies:

\begin{equation}\label{eqn:newderivorder}\displaystyle 0 \leq \min\limits_{j\in\N} \,\mbox{ord}_x \tilde{\mu}_{m_0q, \, p, \, j} \leq m_0q\,\, \nu(\tilde{\psi}_{m_0, \, p}, \, x)\leq nq, \hspace{2ex} x\in \tilde{U}.\end{equation}

\noindent Thus all the common zeroes of the $\tilde{\mu}_{m_0q,p, j}$'s can be removed by deriving them up to order $\geq nq$. With a view to constructing regularising functions, we set the following analogue of the functions $u_m$ defined in (\ref{eqn:regderiv}) for which we obtained a lower bound in Corollary \ref{Cor:below}. While $m=m_0q$, ${\mathcal H}_{\Omega}(m_0q\varphi)$ is now replaced with ${\mathcal H}_{\tilde{U}}(m_0q\tilde{\psi}_{m_0,p})$ having $(\tilde{\mu}_{m_0q,p, j})_{j\in\N}$ as an orthonormal basis. By (\ref{eqn:newderivorder}), $m\nu$ can be replaced with $nq$.

\begin{Def}\label{Def:regabove} For every $m_0, q\in \N$ and $\varepsilon = p\,\frac{n+2}{m_0}$, we let $m=m_0q$ and \\

$\displaystyle \tilde{\psi}_{m_0q, \, p}:=\frac{1}{2m_0q}\, \log\sum\limits_{j=1}^{N_{m_0q, \, p}}\sum\limits_{|\alpha|=0}^{[nq(1+\varepsilon)]}\bigg|\frac{D^{\alpha}\tilde{\mu}_{m_0q,p, j}}{\alpha!}\bigg|^2 + \frac{1}{m_0}\, \log|\tilde{g}_{m_0, p}|.$  \\

\end{Def}

 The advantage over the situation summed up in Corollary \ref{Cor:below} is that the maximal derivation order $[nq(1+\varepsilon)]$ is now small compared with $m_0q$ appearing in the denominator (i.e. $[nq(1+\varepsilon)]/m_0q \rightarrow 0$ when $m_0\rightarrow +\infty$). In view of Lemma \ref{Lem:derivorder}, this is a significant improvement leading to a new regularisation of $\varphi_p$.

\begin{Lem}\label{Lem:localreg} The functions $\psi_{m_0q,\,p}:=\tilde{\psi}_{m_0q, \, p} \circ \mu^{-1}$ defined on $U\setminus V\, {\mathcal I}(m_0(1+\varepsilon)\varphi_p)$ extend to psh functions on $U$ which converge pointwise and in $L^1_{loc}$ topology to $\varphi_p$ when $m_0, q \rightarrow +\infty$. In particular, $dd^c\psi_{m_0q,\,p}$ converges to $dd^c\varphi_p$ as currents.

\end{Lem}

\noindent {\it Proof.} Lemma \ref{Lem:newreg} applied to the orthonormal basis $(\tilde{\mu}_{m_0q,\,p,\,j})_{j\in\N}$ of ${\mathcal H}_{\tilde{U}}(m_0q\, \tilde{\psi}_{m_0,\,p})$ with $m=m_0q$, $\varepsilon_m=\frac{n+2}{m_0}(1+\varepsilon),$ and $C=\frac{n}{n+2}$ gives: 

$\displaystyle \varphi_p \circ \mu(w) - \frac{\log N_{m_0q,\,p}+C_1}{2m_0q}\leq\tilde{\psi}_{m_0q, \, p}(w) \leq$   \\

\noindent\begin{equation}\label{eqn:pestimate}\sup\limits_{|v-w|<2r}\bigg(\psi_{m_0,\,p}(v)+\frac{1}{m_0}\, \log|\tilde{g}_{m_0,\,p}(w)|\bigg) - \frac{[nq(1+\varepsilon)]+n}{m_0q}\, \log r + \frac{1}{m_0q}\, \log C_3,\end{equation}

\noindent for every $w\in \tilde{U}$ and every $0\leq r < \min\{\frac{1}{2}\, \mbox{dist} (w,\, \partial \tilde{U}),\, 1\}$. As every $w\in\tilde{U}\setminus E_{m_0(1+\varepsilon),\,p}$ is the image of a unique $z\in U\setminus V\,{\mathcal I}(m_0q(1+\varepsilon)\varphi_p)$ under $\mu^{-1}$, the conclusion follows. \hfill $\Box$

\vspace{2ex}

We can now estimate the growth of the Monge-Amp\`ere masses of $(dd^c\psi_{m_0q, \, p})^k$, $k=1, \dots , \, n,$ regarded as well-defined smooth forms on their smooth locus $U\setminus V\, {\mathcal I}(m_0(1+\varepsilon)\varphi_p),$ as $m_0, \, q \rightarrow +\infty$.

\begin{Lem}\label{Lem:pmasses} For every $k=1, \dots , n$, the following mass estimate holds: \\

\hspace{6ex} $\displaystyle\int\limits_{U'} (dd^c\psi_{m_0q,\,p})_{ac}^k \wedge \beta^{n-k} \leq C_{m_0}\, (-\log\delta_{m_0})^k,$ \hspace{2ex} $m_0 \gg  1, \, q\in \N$ \\

\noindent on every $U'\Subset U$, with a constant $C_{m_0}>0$ independent of $p, q$. In particular, if $q=q(m_0)$ is chosen so big as to have $\frac{1}{q(m_0)}\, C_{m_0}\, (-\log\delta_{m_0})^n\rightarrow 0$ when $m_0\rightarrow +\infty$, we get: \\

\hspace{6ex} $\displaystyle \lim\limits_{m_0\rightarrow +\infty}\frac{1}{m_0\, q(m_0)}\,\int\limits_{U'} (dd^c\psi_{m_0q,\,p})_{ac}^k \wedge \beta^{n-k}=0$, \hspace{2ex} $k=1, \dots , n,$ \\

\noindent which, in view of the result expected in Theorem \ref{The:globalregularisation} $(c)$ with $m=m_0q$, is satisfactory.

\end{Lem}

\noindent {\it Proof.} As the singularities of $\psi_{m_0q,\,p}$ are analytic, the mass of $(dd^c\psi_{m_0q,\,p})_{ac}^k$ on $U'$ equals the mass of $(dd^c\psi_{m_0q,\,p})^k$ outside the singular locus (i.e. on $U'\setminus V{\mathcal I}(m_0(1+\varepsilon)\varphi_p)$. Let $\tilde{U}'\Subset \tilde{U}$ be the inverse image of $U'$ under $\mu$. If we single out the $C^{\infty}$ part of $\tilde{\psi}_{m_0q,\,p}$ as:

\begin{equation}\label{eqn:smoothpart}\displaystyle\tilde{u}_{m_0q,\,p}:=\frac{1}{2m_0q}\, \log\sum\limits_{j=1}^{N_{m_0q, \, p}}\sum\limits_{|\alpha|=0}^{[nq(1+\varepsilon)]}\bigg|\frac{D^{\alpha}\tilde{\mu}_{m_0q,p, j}}{\alpha!}\bigg|^2,\end{equation}

\noindent we have: \\

$\begin{array}{lll}\displaystyle\int\limits_{U'} (dd^c\psi_{m_0q,\,p})_{ac}^k \wedge \beta^{n-k} & = & \displaystyle\int\limits_{\tilde{U}'\setminus E_{m_0(1+\varepsilon),\,p}} (dd^c\tilde{\psi}_{m_0q,\,p})^k \wedge \mu^{\star}\beta^{n-k} \\

 & = & \displaystyle\int\limits_{\tilde{U}'} (dd^c\tilde{u}_{m_0q,\,p})^k \wedge \mu^{\star}\beta^{n-k}.\end{array}$ \\

\noindent The latter equality above follows from the fact that $dd^c(\frac{1}{m_0}\,\log|\tilde{g}_{m_0,\,p}|)$ equals the current of integration on $\mbox{div}\,\tilde{g}_{m_0,\,p}$ whose support is included in $E_{m_0(1+\varepsilon),\,p}$. As the mass is calculated in the complement of $E_{m_0(1+\varepsilon),\,p}$, $dd^c(\frac{1}{m_0}\,\log|\tilde{g}_{m_0,\,p}|)$ has no contribution. By the Chern-Levine-Nirenberg inequalities ([CLN69]) we get:

\begin{equation}\label{eqn:CLN}\displaystyle\int\limits_{\tilde{U}'} (dd^c\tilde{u}_{m_0q,\,p})^k \wedge \beta^{n-k} \leq C'_{m_0}\, (\sup\limits_{\tilde{U}}|\tilde{u}_{m_0q,\,p}|)^{k}, \hspace{2ex} k=1, \dots, n,\end{equation} 

\noindent with a constant $C'_{m_0}>0$ depending only on $\tilde{U}'\Subset\tilde{U}\Subset \tilde{\Omega}_{m_0}$ (and implicitly on $m_0$). The supremum above is finite since $\tilde{u}_{m_0q,\,p}$ is $C^{\infty}$. Thus controlling the Monge-Amp\`ere masses comes down to controlling the growth of $\sup\limits_{\tilde{U}}|\tilde{u}_{m_0q,\,p}|$ as $m_0, q\rightarrow +\infty$.

\begin{Claim}\label{Claim:avcontrol} There are constants $C''_{m_0},\, A_{m_0}\geq 0$ independent of $p, q$ such that for $q\gg  m_0\gg  p$ we have:\\

 \hspace{6ex} $\displaystyle |\tilde{u}_{m_0q,\,p}|\leq C''_{m_0}\, (-\log\delta_{m_0}) + A_{m_0}$ \hspace{2ex} on $\tilde{U}'\Subset\tilde{U}.$  \\

\end{Claim}

 \noindent As $(\tilde{\mu}_{m_0q,\,p,\,j})_{j\in \N}$ defines an orthonormal basis of ${\mathcal H}_{\tilde{U}}(m_0q\tilde{\psi}_{m_0,\,p})$, we will make use of the results in section \ref{subsection:addefect} to get the desired control. To prove the Claim, we need to find upper and lower bounds for $\tilde{u}_{m_0q,\,p}$. We first settle the upper bound question. We can use the effective version of the subadditivity property of multiplier ideal sheaves involving derivatives. Proposition \ref{Prop:effectiveincl} $(b)$, applied to the orthonormal bases $(\tilde{\mu}_{m_0q,\,p,\,j})_{j\in\N}$ and $(\tilde{\mu}_{m_0,\,p,\,j})_{j\in\N}$ of ${\mathcal H}_{\tilde{U}}(m_0q\tilde{\psi}_{m_0,\,p})$ and ${\mathcal H}_{\tilde{U}}(m_0\tilde{\psi}_{m_0,\,p})$, gives: \\

\noindent$\begin{array}{lll}\displaystyle \sum\limits_{j=0}^{+\infty}|D^{\alpha}\tilde{\mu}_{m_0q,\,j}|^2 & \leq & \displaystyle C_n^{q-1}\, \sum\limits_{j_1, \dots, j_q=0}^{+\infty}|D^{\alpha}(\tilde{\mu}_{m_0,\,j_1}\dots \tilde{\mu}_{m_0,\,j_q})|^2 \\

  & \leq & \displaystyle C_n^{q-1}\, \sum\limits_{j_1, \dots, j_q=0}^{+\infty}\bigg|\sum\limits_{\beta_1 + \dots +\beta_q=\alpha} C^{\alpha}_{\beta_1, \dots , \beta_q} \, D^{\beta_1}\tilde{\mu}_{m_0,\,j_1}\dots D^{\beta_q}\tilde{\mu}_{m_0,\,j_q})\bigg|^2 \\

  & \leq & \displaystyle C_n^{q-1}\, C_{\alpha}\,\bigg(\sum\limits_{j=0}^{+\infty}\sum\limits_{\beta\leq \alpha}|D^{\beta}\tilde{\mu}_{m_0,\,j}|^2\bigg)^q,\end{array}$\\

\noindent for every $\alpha\in \N^n$, where $C_{\alpha}:=\sum\limits_{\beta_1 + \beta_q=\alpha}(C^{\alpha}_{\beta_1, \dots , \beta_q})^2$. This implies that: \\

\hspace{6ex} $\displaystyle \tilde{u}_{m_0q,\,p} \leq \tilde{u}_{m_0,\,p} + \frac{\log N_{m_0,\,p}}{2m_0} + \frac{C_q}{m_0q}$ \hspace{2ex} on $\tilde{U},$ \\

\noindent where $C_q>0$ is a constant depending only on $q$ and $n$ such that $C_q/q$ is bounded as $q\rightarrow +\infty$. This means that the sequence $(\tilde{u}_{m_0q,\,p})_{m_0,\,q\in\N}$ is non-increasing up to constants independent of $q$. It is in particular bounded above if we choose $m_0\gg p$.

 For the much subtler problem of finding a lower bound for $\tilde{u}_{m_0q,\,p}$ (cf. (\ref{eqn:smoothpart})), we use Corollary \ref{Cor:below} to get an orthonormal basis $(\tilde{\mu}_{m_0q,\,p,\,j})_{j\in \N}$ of ${\mathcal H}_{\tilde{U}}(m_0q\tilde{\psi}_{m_0,\,p})$ such that the derivatives up to order $[n(1+\varepsilon)q]$ satisfy: \\

\hspace{6ex}$\displaystyle \tilde{u}_{m_0q,\,p}\geq C_0^{(m_0, p)}\, \log\delta_{m_0}-A_{m_0q}$ \hspace{2ex} on \, $\tilde{U}'\Subset\tilde{U}.$ \\

\noindent As $\varepsilon=p\,\frac{n+2}{m_0}$ (cf. Definition \ref{Def:regabove}), the estimates of Corollary \ref{Cor:below} give, for $q\gg m_0\gg p$: \\

\hspace{6ex} $\displaystyle A_{m_0q}\leq \frac{\log(C_n q^n N_{m_0q,\,p}^2)}{m_0q} + 2(n+2)\frac{p}{m_0}\, \sup\limits_{\tilde{U}'}\tilde{\psi}_{m_0,\,p}\leq A_{m_0},$ \\

\noindent with $C_n>0$ depending on $\tilde{U}',\tilde{U}$ and implicitly on $m_0$. We also have: \\

$\displaystyle 0\leq C_0^{(m_0, p)}\leq C_{\tilde{U}}\, \int_{\tilde{U}}dd^c\tilde{\psi}_{m_0,\,p}\wedge \mu^{\star}\beta^{n-1} \leq C_{\tilde{U}}\, \int_U dd^c\varphi_p\wedge \beta^{n-1}$. \\

\noindent Since $dd^c\varphi_p$ converges weakly to $dd^c\varphi$ on $\Omega$ as $p\rightarrow +\infty$, the sequence $(\int_U dd^c\varphi_p\wedge \beta^{n-1})_{p\in\N}$ is bounded. As $\tilde{U}\Subset \tilde{\Omega}_{m_0}$, $C_{\tilde{U}}$ depends on $m_0$. Thus the constants $C_0^{(m_0, p)}$ are bounded by a constant depending only on $m_0$ and independent of $p$. Moreover, $\delta_{m_0}$ depends a priori on $p$ as the minimum mutual distance of points where $\varphi_p\circ \mu$ has large enough Lelong numbers (see the proof of Proposition \ref{Prop:derivlowerbound} for the definition of $\delta_{m_0}$), but it can actually be taken independent of $p$ since the singularities of all $\varphi_p$'s are among those of $\varphi$ and the Lelong numbers of $\varphi_p$ are within $n/p$ of those of $\varphi$ (see (\ref{eqn:Lelongnumbers})). 

 This completes the proof of Claim \ref{Claim:avcontrol}. Lemma \ref{Lem:pmasses} follows from estimate (\ref{eqn:CLN}) and Claim \ref{Claim:avcontrol} by taking $C_{m_0}=C'_{m_0}\, C^{''k}_{m_0}$ and absorbing $A_{m_0}$ in $C_{m_0}$. \hfill $\Box$

\vspace{2ex}

  Estimate (\ref{eqn:pestimate}) satisfied by $\tilde{\psi}_{m_0q,\, p}$ (cf. the proof of Lemma \ref{Lem:localreg}), combined with estimate (\ref{eqn:approximation}) satisfied by $\varphi_p$, implies that $\psi_{m_0q,\,p}$ converges pointwise and in $L^1_{loc}$ topology to $\varphi$ as $m_0,\, q,\, p\rightarrow +\infty$. Now Lemma \ref{Lem:localreg} defines a function $\psi_{m_0q,p}$ on each open set of the covering $(U_l)_{1\leq l \leq N}$ of $\Omega$ (take $U=U_l$). Taking $m=m_0q$ with $q\gg m_0\gg p$, we can patch together the psh functions:

\hspace{25ex}$\psi_{m,\,p}$ \\

\noindent defined on various open subsets $U_l\subset\Omega$, ${1\leq l \leq N},$ into regularising functions for $\varphi$ on $\Omega$. The following is a local version of Theorem \ref{The:globalregularisation}.

\begin{Prop}\label{Prop:localapsh} Let $\varphi$ be an arbitrary psh function on a bounded pseudoconvex domain $\Omega\Subset \C^n$. There exists a sequence of almost psh functions $(\psi_m)_{m\in\N}$ with analytic singularities such that $\psi_m$ converges pointwise on $\Omega$ and in $L^1_{loc}$ topology to $\varphi$ as $m\rightarrow +\infty$, each $\psi_m$ is smooth on $\Omega\setminus V{\mathcal I}(m\varphi)$, and the following hold: \\

\noindent $(a)$\, $dd^c\psi_m \geq -\frac{C}{m}\, \beta$ \,\, for some constant $C>0$ independent of $m$; \\

\noindent $(b)$\, $\nu(\varphi,\, x)-\varepsilon_m\leq \nu(\psi_m,\,x)\leq \nu(\varphi,\, x)$, \hspace{1ex} $x\in\Omega,\, m\in\N$, for some $\varepsilon_m\downarrow 0$;\\

\noindent $(c)$\, $\displaystyle\lim\limits_{m\rightarrow +\infty}\frac{1}{m}\, \int\limits_B(dd^c\psi_m+\frac{C}{m} \, \beta)_{ac}^k\wedge \beta^{n-k}=0,$ \hspace{2ex} $k=1,\dots , n,$ \\

\noindent for every relatively compact open subset $B\Subset \Omega$.

\end{Prop}

\noindent {Proof.} For any $x_0\in\Omega$, adding $|z-x_0|^2$ to $\varphi$ if necessary, we can achieve that $i\partial\bar{\partial}\varphi \geq \beta$ and thus the results of this section apply. Let $(\varphi_p)_{p\in\N}$ be the Demailly regularisation (cf. (\ref{eqn:Demaillyapprox})) of $\varphi$. For every $l=1, \dots , N$, we define $\psi_{m_0q,\,p}^{(l)}$ on $U_l$ like in Definition \ref{Def:regabove} and Lemma \ref{Lem:localreg} with $U$ in place of $U_l$. We can consider finite covers $(U'_l)_{1\leq l \leq N}$, $(U''_l)_{1\leq l \leq N}$ and $(U_l)_{1\leq l \leq N}$ of $\Omega$ by concentric balls of radii $\delta$, $\frac{3}{2}\delta$ and respectively $2\delta$. We then use the (essentially well-known) patching procedure recalled in section \ref{subsection:global} to set:\\

\hspace{6ex}$\psi_{m_0q,\,p}(z):=\sup\limits_{U''_l\ni z}\bigg(\psi_{m_0q,\,p}^{(l)}(z) + \frac{C(\delta)}{m_0q}(\delta^2-|z^{(l)}|^2)\bigg),$ \\

\noindent and show, by means of H\"ormander's $L^2$ estimates, that there is a constant $C(\delta)>0$ depending only on $\delta$ for which the patching condition (\ref{eqn:patching}) holds. Here $z^{(l)}$ is a coordinate system centred at the centre of $U^{(l)}$ and H\"ormander's $L^2$ estimates are applied with the analogous weight to (\ref{eqn:weight}) in which $\log|z-x_0|$ has now the coefficient $n+[n(1+\varepsilon)q]$ to force the solution of the equation (\ref{eqn:horeqn}) to vanish at $x_0$ to order $[n(1+\varepsilon)q]$. For every $m_0$, we choose $q=q(m_0)\gg m_0$ like in Lemma \ref{Lem:pmasses}. For $m=m_0q$, we choose an increasing sequence of integers $p_m\rightarrow +\infty$ with $p_m<<m_0$ and set: \\

\hspace{25ex} $\psi_m:=\psi_{m,\,p_m},$ \\

\noindent which is easily seen to satisfy the requirements thanks to the results obtained in the present section.  \hfill $\Box$

\vspace{2ex}

\noindent {\it End of Proof of Theorem \ref{The:globalregularisation}.}  Theorem \ref{The:globalregularisation} stated in the Introduction now follows by patching together functions analogous to $\psi_m$ obtained in the above Proposition \ref{Prop:localapsh} on various open sets contained in coordinate patches covering $X$. The patching procedure, recalled in section \ref{subsection:global} and applied there and in the above proof, can be repeated. A word of explanation is in order here. This patching procedure, based on H\"ormander's $L^2$ estimates, was used in [Dem92] to patch together regularising functions of the type (\ref{eqn:Demaillyapprox}) when no derivatives are involved (see Proposition 3.7. in [Dem92]). In order to obtain smooth regularisations, derivatives (or jets) had to be introduced and a new patching procedure, based on Skoda's $L^2$ estimates, was devised in [Dem92, section 5] to handle this situation. It is worth stressing that this latter, more powerful patching procedure was necessary there to handle derivatives of order up to $[cm]$ in the definition of the $m^{th}$ regularising function. In our present case, we only derive up to order $[C\varepsilon_m m]$ with $\varepsilon_m\downarrow 0$ (see Lemma \ref{Lem:newreg} and Definition \ref{Def:regabove}) and the former patching procedure can be used with minor modifications as explained in section \ref{subsection:global} when the derivation order was the constant $1$.   \hfill $\Box$

\section{Singular hermitian metrics and big line bundles}\label{section:bigness}

 We are now in a position to prove the analytic characterisation of the volume of a line bundle spelt out in Theorem \ref{The:bigness} as a geometric application of our current regularisation Theorem \ref{The:globalregularisation} with controlled Monge-Amp\`ere masses. As mentioned in the Introduction, the case of a non-K\"ahler compact manifold $X$ is new. 

 We first briefly review the set-up (cf. [Dem85], [Bon98], [Bou02]). Let $(L, \, h)$ be a holomorphic line bundle over a compact Hermitian manifold $(X, \, \omega)$ equipped with a possibly singular Hermitian metric $h$. Let $T:=i\Theta_h(L)$ be the curvature current associated with $h$. There is a global representation of $T$ as $T=\alpha + dd^c\varphi$ with a global $C^{\infty}$ $(1, \, 1)$-form $\alpha$ on $X$. For every $q=0, 1, \dots ,\, n$, define the $q$-index set of $T$ as the open subset $X(q, \, T)$ of $X$ consisting of those points $x$ such that $T_{ac}(x)$ has precisely $q$ negative and $n-q$ positive eigenvalues. Let $X(\leq q, \, T):= X(0, \, T) \cup \dots \cup  X(q, \, T)$. For every $m\in \N^{\star}$, consider the singular metric $h^m$ on $L^m$ induced by $h$. This means that if $h=e^{-\varphi}$ on an open subset $U\subset X$ on which $L$ is trivial, $h^m$ is defined as $h^m=e^{-m\varphi}$ on $U$. Now suppose that $T:=i\Theta_h(L) \geq -C\, \omega$ for some constant $C>0$ (i.e. $T$ is almost positive and $\varphi$ is almost psh). Then the associated multiplier ideal sheaf ${\mathcal I}(h^m)$ is the coherent subsheaf of ${\mathcal O}_X$ defined locally as ${\mathcal I}(h^m)_{|U}={\mathcal I}(m\varphi)$ (see (\ref{eqn:misdef})).

 Demailly's holomorphic Morse inequalities ([Dem85]) for smooth metrics $h$ were generalised by Bonavero ([Bon98]) to the case of singular metrics $h$ with analytic singularities in the form of the following asymptotical estimates for the cohomology group dimensions of the twisted coherent sheaves ${\mathcal O}_X(L^m)\otimes {\mathcal I}(h^m)$:  \\

$\displaystyle\sum\limits_{j=0}^q (-1)^{q-j}\, h^j(X, \, {\mathcal O}_X(L^m)\otimes {\mathcal I}(h^m)) \leq \frac{m^n}{n!} \, \int\limits_{X(\leq q, \, T)}(-1)^q \, T_{ac}^n + o(m^n)$, \\

\noindent as $m\rightarrow \infty$, for all $q=1, \dots , \, n$. For $q=1$, we get: \\

\noindent $\displaystyle h^0(X, \, {\mathcal O}_X(L^m)\otimes {\mathcal I}(h^m)) - h^1(X, \, {\mathcal O}_X(L^m)\otimes {\mathcal I}(h^m)) \geq \frac{m^n}{n!}\, \int\limits_{X(\leq 1, \, T)} T_{ac}^n + o(m^n)$. \\

\noindent As \hspace{2ex} $h^0(X, \, {\mathcal O}_X(L^m))  \geq  h^0(X, \, {\mathcal O}_X(L^m)\otimes {\mathcal I}(h^m))$   \\

\hspace{16ex}  $ \geq  h^0(X, \, {\mathcal O}_X(L^m)\otimes {\mathcal I}(h^m))- h^1(X, \, {\mathcal O}_X(L^m)\otimes {\mathcal I}(h^m)),$ \\

\noindent we infer the following lower bound for the volume of $L$: 

\begin{equation}\label{eqn:volbound}\displaystyle v(L)\geq \int\limits_{X(\leq 1, \, T)} T_{ac}^n,\end{equation}

\noindent for every almost positive closed current $T$ with analytic singularities (if any) in $c_1(L)$. \\

\noindent {\it Proof of Theorem \ref{The:bigness}.} Clearly, it is enough to prove the equality characterising the volume, as the bigness criterion is an immediate consequence of it. The inequality ``$\leq$'' bounding the volume above can be proved as in [Bou02] since $X$ is Moishezon when $v(L)>0$ and can be modified into a projective manifold. If $v(L)=0$, the inequality ``$\leq$'' is obvious.

 Thus proving Theorem \ref{The:bigness} boils down to obtaining a lower bound for the volume of $L$ in terms of curvature currents. In the light of the above explanations, this can be seen as singular Morse inequalities for arbitrary singularities. Let $T:=i\Theta_h(L)\geq 0$ be the curvature current associated with a singular Hermitian metric $h$ with arbitrary singularities on $L$. If no positive current exists in $c_1(L)$, there is nothing to prove. By Theorem \ref{The:globalregularisation}, there exist regularising currents with analytic singularities $T_m\rightarrow T$ in $c_1(L)$ such that $T_m \geq - \frac{C}{m}\, \omega$ for some constant $C>0$ independent of $m$ which satisfy the Monge-Amp\`ere mass condition $(c)$. Furthermore, Theorem 2.4 in [Bou02, p. 1050] asserts that we can modify our sequence $(T_m)_{m\in \N}$ in such a way that besides all its properties it also satisfies:

 \begin{equation}\label{eqn:ac}T_m(x) \rightarrow T_{ac}(x) \hspace{2ex} \mbox{as} \hspace{2ex} m\rightarrow +\infty, \hspace{2ex} \mbox{for almost every} \hspace{2ex} x\in X.\end{equation}

\noindent A word of explanation is in order here. Property (\ref{eqn:ac}) is achieved in [Bou02] by combining the potentials $\psi_m$ of the currents $T_m$ with those of another regularising sequence $\alpha_m$ of smooth forms constructed in [Dem82]. For tubular neighbourhoods $W_m\Subset U_m$ of the singular locus of $\psi_m$, and for arbitrary sequences $\delta_m\downarrow 0$ and $C_m\uparrow +\infty$, $\psi_m$ is replaced with $(1-\delta_m)\psi_m-C_m$ on $U_m$, and with a regularised maximum function of $(1-\delta_m)\psi_m-C_m$ and the smooth potential of $\alpha_m$ on $X\setminus W_m$. In the proof of Theorem \ref{The:globalregularisation} above, we have constructed the potentials $\psi_m$ of $T_m$ to be ``not too small'' (after removing their singularities on a modification) and hence their moduli, which control the Monge-Amp\`ere masses by the Chern-Levine-Nirenberg inequality, to be ``not too large''. It is clear that this property is preserved by the regularised maximum construction of [Bou02] and hence so is the Monge-Amp\`ere mass control of Theorem \ref{The:globalregularisation} (c).

\noindent As explained above, by the Morse inequalities applied to $L$ with $T_m\in c_1(L)$ as curvature current with analytic singularities, we get (cf. (\ref{eqn:volbound})): \\

 $\displaystyle v(L) \geq \int\limits_{X(\leq 1, \, T_m)} T_{m, \, ac}^n = \int\limits_{X(0, \, T_m)} T_{m, \, ac}^n + \int\limits_{X(1, \, T_m)} T_{m, \, ac}^n$ \hspace{2ex} for every $m\in \N.$ \\

\noindent On the other hand, the proof of Proposition 3.1. in [Bou02, p. 1052-53] uses the Fatou lemma to derive the following inequality from property (\ref{eqn:ac}): \\

\hspace{6ex} $\displaystyle \liminf\limits_{m\rightarrow +\infty} \int_{X(0, \, T_m)}T_{m, \, ac}^n \geq \int_{X(T, \, 0)}T_{ac}^n = \int_XT_{ac}^n$. \\

 \noindent Thus, to prove the Morse-type inequality ``$\geq$'' it is enough to show that $\lim\limits_{m\rightarrow +\infty}\int\limits_{X(1, \, T_m)} T_{m, \, ac}^n = 0$. Note that on the open set $X(1, \, T_m)$ we have: \\

\hspace{6ex} $\displaystyle 0 \leq - T_{m, \, ac}^n \leq \displaystyle n \, \frac{C}{m} \, (T_{m, \, ac} + \frac{C}{m}\, \omega)^{n-1} \wedge \omega$. \\

\noindent  It is thus enough to show that \\

\hspace{6ex} $\displaystyle \lim\limits_{m\rightarrow +\infty} \frac{C}{m} \,\int\limits_X  (T_{m, \, ac} + \frac{C}{m}\, \omega)^{n-1}\wedge\omega =0$. \\

\noindent Since $\int\limits_X  (T_{m, \, ac} + \frac{C}{m}\, \omega)^{n-1}\wedge\omega = \int\limits_{X\setminus V{\mathcal I}(mT)} (T_m + \frac{C}{m}\, \omega)^{n-1}\wedge\omega$, this is precisely the Monge-Amp\`ere mass property obtained in Theorem \ref{The:globalregularisation} $(c)$. The proof is thus the same as in the K\"ahler case settled in [Bou02] once we have obtained Theorem \ref{The:globalregularisation} which is new in the non-K\"ahler context.    \hfill $\Box$

\section{Appendix: auxiliary results}\label{subsection:appendix}

For the sake of perspicuity and completeness, we spell out in this Appendix two results related to the $L^2$ estimates that were used in the foregoing sections. In section \ref{subsection:addefect}, an essential use was made of Skoda's $L^2$ division theorem in the following form.

\begin{The}(Skoda [Sko72b])\label{The:Skoda} Let $\varphi$ be a psh function on a pseudoconvex open set $\Omega\subset \C^n$, and let $g_1, \dots , g_N$ be (possibly infinitely many) holomorphic functions on $\Omega$. Set $r:=\min\{N-1, \, n\}$ and $|g|^2=\sum\limits_{j=1}^N|g_j|^2$. Then, for every holomorphic function $f$ on $\Omega$ satisfying:

\begin{equation}\nonumber\int\limits_{\Omega}|f|^2\, |g|^{-2(r+s+\alpha)}\, e^{-2\varphi}\, dV_n <+\infty, \hspace{3ex} \alpha>0, \,\, s\in\N^{\star},\end{equation}

\noindent there exist holomorphic functions $h_L$ on $\Omega$ for all $L=(l_1, \dots , l_s)\in \{1, \dots , N\}^s$ such that:

\begin{equation}\nonumber f=\sum\limits_Lh_L\, g^L \hspace{2ex} \mbox{on}\,\, \Omega, \hspace{5ex} \mbox{with} \hspace{2ex} g^L=g_{l_1}\dots g_{l_s},\end{equation}

\begin{equation}\nonumber \int\limits_{\Omega}\sum\limits_L|h_L|^2\, |g|^{-2(r+\alpha)}\, e^{-2\varphi}\, dV_n \leq \frac{\alpha +s}{\alpha}\, \int\limits_{\Omega}|f|^2\, |g|^{-2(r+s+\alpha)}\, e^{-2\varphi}\, dV_n.\end{equation}

\end{The}

\noindent This statement, in which the original function $f$ is divided by products of $s$ functions $g_l$ if it satisfies an appropriate $L^2$ condition depending on $s\in\N^{\star}$, follows straightforwardly by induction on $s$ from Skoda's original statement given in [Sko72b] for $s=1$. (See e.g. [Dem92, Corollary A.5., p. 407] for the present statement).

\vspace{2ex}

 The other ingredient listed in this Appendix is the proof of Lemma \ref{Lem:locglobal} used in section \ref{subsection:reg}. \\

\noindent {\bf Proof of Lemma \ref{Lem:locglobal}}  The restriction to $B$ clearly defines an injection ${\mathcal H}_{\Omega}(m\varphi)\hookrightarrow {\mathcal H}_{B}(m\varphi)$ and the image of the unit ball of the first space is contained in the unit ball of the second space. Thus the first inequality holds on $B$. For the second inequality, let $x\in B_0$ be a fixed point,  and suppose, for instance, that $B_0=B(x, \, r/2)$. Let $\theta$ be a $C^{\infty}$ cut-off function such that \\

 $\mbox{Supp}\, \theta \Subset B(x, \, r) \Subset B,$ \hspace{2ex} $\theta\equiv 1$ on $B(x, \, r/2),$ \hspace{2ex} $0\leq \theta \leq 1$ on $\C^n.$ \\

\noindent Let $f\in {\mathcal O}(B)$ such that $\int_B |f|^2\, e^{-2m\varphi}=1$ be an arbitrary element in the unit sphere of ${\mathcal H}_B(m\varphi)$. We need produce a global holomorphic function on $\Omega$ all of whose derivatives up to order $p$ assume the same values at $x$ as those of $f$. We can use H\"ormander's $L^2$ estimates (cf. [Hor65]) to solve the equation: \\

\hspace{15ex}$\bar{\partial}u = \bar{\partial}(\theta\, f)$ \hspace{2ex} on $\Omega$ \\

 \noindent with the strictly psh weight $m\varphi + (n+p)\log|z-x| + |z-x|^2$. There exists a $C^{\infty}$ solution $u$ satisfying the estimate: \\

\hspace{6ex} $\displaystyle \int\limits_{\Omega}\frac{|u|^2}{|z-x|^{2(n+p)}}\, e^{-2m\varphi}\, e^{-2|z-x|^2} \leq 2 \int\limits_{\Omega}\frac{|\bar{\partial}\theta|^2\, |f|^2}{|z-x|^{2(n+p)}}\, e^{-2m\varphi}\, e^{-2|z-x|^2}.$ \\

\noindent The non-integrability of $|z-x|^{-2(n+k)}$, $k=0,1, \dots , p$, near $x$ implies that $D^{\alpha}u(x)=0$, $0\leq |\alpha|\leq p$. If we set: \\

\hspace{6ex} $F:=\theta \, f -u\in {\mathcal O}(\Omega),$ \\

\noindent we have $D^{\alpha}F(x)=D^{\alpha}f(x)$, $0\leq |\alpha|\leq p,$ and \\

\hspace{3ex} $\displaystyle \int\limits_{\Omega}|F|^2\, e^{-2m\varphi} \leq 2\bigg(1+C_n\, \frac{d^{2n}\, e^{2d^2}}{r^{2(n+1)}}\bigg)\,(d/r)^p\, \int\limits_B |f|^2\, e^{-2m\varphi}=C_{n, \, d,\, r}(d/r)^p,$ \\

\noindent since $\int_B|f|^2\, e^{-2m\varphi}=1$, with a constant $C_n>0$ depeding only on $n$ and $C_{n, \, d, \, r}$ denoting the double of the parenthesis above. This means that $F/\sqrt{C_{n, \, d, \, r}(d/r)^p}$ belongs to the unit ball of ${\mathcal H}_{\Omega}(m\varphi)$. The second inequality follows from the following expressions holding at every $x\in B$: \\

$\displaystyle B^{(p)}_{m\varphi, \, \Omega}(x)=\sum\limits_{|\alpha|=0}^{p}\sup\limits_{g\in \bar{B}_{m, \, \Omega}(1)}|D^{\alpha}g(x)|^2$, \hspace{2ex} $\displaystyle B^{(p)}_{m\varphi, \, B}(x)=\sum\limits_{|\alpha|=0}^{p}\sup\limits_{f\in \bar{B}_{m, \, B}(1)}|D^{\alpha}f(x)|^2,$ \\

\noindent where $\bar{B}_{m, \, \Omega}(1)$ and $\bar{B}_{m, \, B}(1)$ are the closed unit balls of ${\mathcal H}_{\Omega}(m\varphi)$ and respectively ${\mathcal H}_B(m\varphi)$.  The calculation details are left to the reader. \hfill $\Box$

\vspace{3ex}

\noindent{\bf References.}

\vspace{2ex}

\noindent [Bou02] \, S. Boucksom --- {\it On the Volume of a Line Bundle} --- Internat. J. of Math. {\bf 13} (2002), no. 10, 1043-1063, MR1945706, Zbl 1101.14008.

\vspace{1ex}

\noindent [Bon95] \, L. Bonavero --- {\it  In\'egalit\'es de Morse et vari\'et\'es de Moishezon} ---  Th\`ese de l'Universit\'e Joseph Fourier, Grenoble (1995), available at http://www-fourier.ujf-grenoble.fr/~bonavero/ or as arXiv preprint alg-geom/9512013.

\vspace{1ex}

\noindent [Bon98] \, L. Bonavero --- {\it In\'egalit\'es de Morse holomorphes singuli\`eres} --- J. Geom. Anal.  {\bf 8} (1998), 409-425, MR1707735, Zbl 0966.32011.

\vspace{1ex}

\noindent [CLN69] \, S.S. Chern, H.I. Levine, L. Nirenberg --- {\it Intrinsic Norms on a Complex Manifold} --- Global Analysis (papers in honour of K. Kodaira), p. 119-139, Univ. of Tokyo Press, Tokyo, 1969, MR0254877, Zbl 0202.11603.

\vspace{1ex}

\noindent [Dem82] \, J.- P. Demailly --- {\it Estimations $L^2$ pour l'op\'erateur $\bar{\partial}$ d'un fibr\'e vectoriel holomorphe semi-positif au dessus d'une vari\'et\'e k\"ahlerienne compl\`ete} --- Ann. Sci \'Ecole Norm. Sup. {\bf 15} (1982), 457-511, MR0690650, Zbl 0507.32021.

\vspace{1ex}

\noindent [Dem85] \, J.- P. Demailly --- {\it Champs magn\'etiques et in\'egalit\'es de Morse pour la $d''$-cohomologie} --- Ann. Inst. Fourier (Grenoble) {\bf 35} (1985), 189-229, MR0812325, Zbl 0595.58014.

\vspace{1ex}

\noindent [Dem92] \, J. -P. Demailly --- {\it Regularization of Closed Positive Currents and Intersection Theory} --- J. Alg. Geom., {\bf 1} (1992), 361-409, MR1158622, Zbl 0777.32016.

\vspace{1ex}

\noindent [Dem93] \, J.- P. Demailly --- {\it A Numerical Criterion for Very Ample Line Bundles} --- J. Differential Geom.  {\bf 37} (1993), 323-374, MR1205448, Zbl 0783.32013.

\vspace{1ex}

\noindent [Dem 97] \, J.-P. Demailly --- {\it Complex Analytic and Algebraic Geometry}---http://www-fourier.ujf-grenoble.fr/~demailly/books.html



\vspace{1ex}

\noindent [DEL00] \, J.-P. Demailly, L. Ein, R. Lazarsfeld --- {\it A Subadditivity Property of Multiplier Ideals}--- Michigan Math. J., {\bf 48} (2000), 137-156, MR1786484, Zbl 1077.14516.

\vspace{1ex}

\noindent [DPS01] \, J. -P. Demailly, T. Peternell, M. Schneider --- {\it Pseudo-effective Line Bundles on Compact K\"ahler Manifolds} --- Internat. J. Math. {\bf 12} (2001), no. 6, 689--741, MR1875649, Zbl 1111.32302.

\vspace{1ex}

\noindent [DP04] \, J.-P. Demailly, M. Paun --- {\it Numerical Characterization of the K\"ahler Cone of a Compact K\"ahler Manifold} --- Ann. of Math. (2) {\bf 159}  (2004),  no. 3, 1247--1274, MR2113021, Zbl 1064.32019.

\vspace{1ex}

\noindent [Dra01] \, D. Drasin ---{\it Approximation of Subharmonic Functions with Applications}--- Approximation, complex analysis, and potential theory (Montreal, QC, 2000), 163--189, NATO Sci. Ser. II Math. Phys. Chem., 37, Kluwer Acad. Publ., Dordrecht, 2001, MR1873588, Zbl 0991.31001.

\vspace{1ex}

\noindent [GR70] \, H. Grauert, O. Riemenschneider --- {\it Verschwindungss\"atze f\"ur analytische Kohomologiegruppen auf komplexen R\"aumen} --- Invent. Math. {\bf 11} (1970), 263-292, MR0302938, Zbl 0202.07602.

\vspace{1ex}

\noindent [Hor65] \, L. H\"ormander --- {\it $L^2$ Estimates and Existence Theorems for the $\bar{\partial}$ Operator} --- Acta Math. {\bf 113} (1965) 89-152, MR0179443, Zbl 0158.11002.

\vspace{1ex}

\noindent [JS93] \, S. Ji, B. Shiffman --- {\it Properties of Compact Complex manifolds Carrying Closed Positive Curents} --- J.Geom. Anal. {\bf 3}, No. 1, (1993), 37-61, MR1197016, Zbl 0784.32009.

\vspace{3ex}

\noindent [Kis84] C. O. Kiselman --- {\it Sur la d\'efinition de l'op\'erateur de Monge-Amp\`ere} --- Lecture Notes in Math. {\bf 1094}, Springer Verlag (1984), 139-150, MR0773106, Zbl 0562.35021.

\vspace{1ex}

\noindent [Nad90])  \, A. M. Nadel ---{\it Multiplier Ideal Sheaves and Existence of K\"ahler-Einstein Metrics of Positive Scalar Curvature}---Proc. Nat. Acad. Sci. U.S.A.  {\bf 86} (1989)  7299-7300, and Ann. of Math.  {\bf 132}  (1990) 549-596, MR1078269, Zbl 0731.53063.

\vspace{1ex}

\noindent [OT87] \, T. Ohsawa, K. Takegoshi --- {\it On The Extension of $L^2$ Holomorphic Functions}--- Math. Zeitschrift {\bf 195} (1987) 197-204, MR0892051, Zbl 0625.32011.

\vspace{1ex}

\noindent [Ohs88] \, T. Ohsawa --- {\it On the Extension of $L^2$ Holomorphic Functions, II}---Publ. RIMS, Kyoto Univ. {\bf 24} (1988) 265-275, MR0944862, Zbl 0653.32012.

\vspace{1ex}

\noindent [Pop04] \, D. Popovici --- {\it Estimation effective de la perte de positivit\'e dans la r\'egularisation des courants} --- C. R. Acad. Sci. Paris, Ser. I 338 (2004) 59-64, MR2038086, Zbl 1063.32002.

\vspace{1ex}

\noindent [Pop06] \, D. Popovici --- {\it Effective Local Finite Generation of Multiplier Ideal Sheaves} --- arXiv e-print math.CV/0603734.

\vspace{1ex}

\noindent [Siu74] \, Y. T. Siu --- {\it Analyticity of Sets Associated to Lelong Numbers and the Extension of Closed Positive Currents} --- Invent. Math. , {\bf 27} (1974), 53-156, MR0352516, Zbl 0289.32003.

\vspace{1ex}

\noindent [Siu75] \, Y. T. Siu --- {\it Extension of Meromorphic Maps into K\"ahler Manifolds}--- Ann. of Math. (2)  {\bf 102}  (1975), no. 3, 421--462, MR0463498, Zbl 0318.32007.

\vspace{1ex}

\noindent   [Siu85] \, Y. T. Siu --- {\it Some Recent Results in Complex Manifolds Theory Related to Vanishing Theorems for the Semipositive Case} --- L. N. M.  {\bf 1111}, Springer-Verlag, Berlin and New-York, (1985), 169-192, MR0797421, Zbl 0577.32032.

\vspace{1ex}

\noindent [Sko72a]  \, H. Skoda---{\it Sous-ensembles analytiques d'ordre fini ou infini dans $\C^n$}--- Bull. Soc. Math. France  {\bf 100}  (1972)  353-408, MR0352517, Zbl 0246.32009.

\vspace{1ex}

\noindent [Sko72b]  \, H. Skoda---{\it Application des techniques $L^2$ \`a la th\'eorie des id\'eaux d'une alg\`ebre de fonctions holomorphes avec poids} --- Ann. Sci. \'Ecole Norm. Sup. (4) 5 (1972), 545-579, MR0333246, Zbl 0254.32017.

\vspace{1ex}

\noindent [Yul85] \, R. S. Yulmukhametov ---{\it Approximation of Subharmonic Functions}---Anal. Math.  {\bf 11} (1985), no. 3, 257-282 (in Russian), MR0822590, Zbl 0594.31005.

\vspace{3ex}

\noindent Dan Popovici

\noindent Universit\'e Paul Sabatier, Laboratoire \'Emile Picard, 118 route de Narbonne, 31062, Toulouse Cedex 4, France.

\noindent E-mail: popovici@math.ups-tlse.fr

\end{document}